\documentclass[11pt,leqno]{amsart}
\usepackage{amsmath,amsfonts,latexsym,graphicx,amssymb,url, color}
\usepackage{txfonts} 

\usepackage[hmargin=2.54 cm,vmargin=2.54 cm]{geometry}

\usepackage{color} 

\newcommand{\R}{{\mathbb R}} 

\newcommand{\U}{{\mathbb U}}

\newcommand{\K}{{\mathbb K}}

\newcommand{\N}{{\mathbb N}}

\newcommand{\e}{\varepsilon} 

\newcommand{\A}{\mathbb{A}}

\newcommand\norm[1]{\left\| #1\right\|}

\newcommand{\M}{{\mathcal M}}

\newcommand{\wei}[1]{\langle #1 \rangle}

\newcommand{\F}{{\mathbf F}}

\newcommand{\bu}{\mathbf{u}}

\newcommand{\tA}{{\tilde{\A}}}

\newtheorem{theorem}{Theorem}[section]
\newtheorem{definition}{Definition}[section]
\newtheorem{remark}{Remark}
\newtheorem{lemma}[theorem]{Lemma}
\newtheorem{proposition}[theorem]{Proposition}

\numberwithin{equation}{section}

\newcommand{\beq}{\begin{equation}}
\newcommand{\eeq}{\end{equation}}

\definecolor{darkred}{rgb}{.70,.12,.20}

\definecolor{darkgreen}{rgb}{.20,.52,.14}

\title[Gradient estimates, Singular quasi-linear parabolic equations] {Regularity gradient estimates for weak solutions of singular quasi-linear parabolic equations}
\author{Tuoc  Phan}
\address{Department of Mathematics, University of Tennessee, Knoxville, 227 Ayres Hall, 1403 Circle Drive, Knoxville, TN 37996, U.S.A.}
\email{phan@math.utk.edu}
\begin{document}
\begin{abstract} This paper studies the Sobolev regularity estimates of weak solutions of a class of singular quasi-linear elliptic problems of the form $u_t - \textup{div}[\A(x, t,u, \nabla u)]= \textup{div}[\F]$ with homogeneous Dirichlet boundary conditions over bounded spatial domains. Our main focus is on the case that the vector coefficients $\A$ are discontinuous and singular in $(x,t)$-variables, and dependent on the solution $u$. Global and interior weighted $W^{1,p}(\Omega, \omega)$-regularity estimates are established for weak solutions of these equations, where $\omega$ is a weight function in some Muckenhoupt class of weights. The results obtained are even new for linear equations, and for $\omega =1$, because of the singularity of the coefficients in $(x,t)$-variables.
\end{abstract}

\maketitle
\section{Introduction and main results} \label{Intro-sec}
Let  $\Omega$ be a bounded domain in $\R^n$ with boundary $\partial \Omega$, let $T >0$, and $\mathbb{K}$ be an open interval in $\mathbb{R}$ including $\K =\R$. The theme of the paper is to establish the interior and global weighted $W^{1,p}$-regularity estimates for weak solutions of the homogeneous Dirichlet boundary value problem  with singular coefficients
\begin{equation}  \label{main-eqn}
\left\{
\begin{array}{cccl}
u_t - \text{div} [\A(x, t,u, \nabla u)] & = & \text{div}[{\bf F}(x,t )]  &\quad \text{in} \quad \Omega_T := \Omega \times (0,T),\\
  u(x,t) & = & 0 & \quad \text{on} \quad \partial  \Omega \times (0,T), \\
  u(x,0) & = & u_0(x) & \quad \text{for} \quad  x \in \Omega.
  \end{array}
\right.
\end{equation}
where ${\bf F}:\Omega_T \rightarrow \mathbb{R}^{n}$ is a given measurable vector field, $u_0: \Omega \rightarrow \mathbb{R}$ is some given measurable function, and $\A : \Omega_T \times \K \times \R^n \rightarrow \R^n$ is measurable and satisfies some natural growth assumptions. 
 
Our main interest in this paper is to establish the $W^{1,p}$-regularity estimates of Calder\'on-Zygmund type for weak solutions of the nonlinear equations \eqref{main-eqn} when the coefficient $\A$ is singular in $(x,t)$-variable and dependent on the solution $u$. As we will explain shortly, treating the singularity of $\A$ in our paper requires some new and very nice ingredients from harmonic analysis theory. Moreover, the dependent of the coefficient $\A$ on $u$ also creates some serious obstacles due to the homogeneity in Calder\'on-Zygmund type estimates which is not available in \eqref{main-eqn}. 

We now precisely state our main assumptions on the coefficients $\A$. We assume that the vector field function $\A: \Omega_T \times \mathbb{K} \times  \mathbb{R}^{n}~\rightarrow~\mathbb{R}^{n}$ is  a Carath\'{e}odory mapping satisfying the following natural  growth conditions
\begin{align} \label{Caratho-1}
& \A(x, t, s, \cdot) \quad \text{is continuous  on} \ \mathbb{R}^n, \ \text{for almost everywhere} \ (x,t, s) \in \Omega_T \times \K, \\
\label{Caratho-2}
& \A(\cdot, \cdot, \cdot, \eta) \quad \text{is measurable for each fixed} \ \eta \in  \mathbb{R}^n, \\
\label{up-ellip}
& \Lambda^{-1} |\xi|^2 \leq \wei{\A(x,t, s, \xi), \xi},  \ \text{for almost everywhere} \ (x,t,s) \in \Omega_T \times \K, \ \text{for all} \   \xi \in  \mathbb{R}^n,
\end{align}
where $\Lambda >0$ is a fixed constant. We assume further that the vector field $\A$ is asymptotically Uhlenbeck in the following sense: There is a measurable matrix $\tilde{\A}: \Omega_T \rightarrow \mathbb{R}^{n \times n}$, and a bounded continuous function $\omega_0 : \overline{\mathbb{K}} \times [0, \infty) \rightarrow [0, \infty)$ such that
\begin{equation} \label{asymp-Uh}
|\A(x,t, s, \eta) - \tilde{\A}(x,t) \eta| \leq \omega_0 (s, |\eta|)\Big[1+ |\eta| \Big], \quad \text{for almost everywhere} \quad (x,t,s) \in \Omega_T \times \K, \quad \forall \  \eta \in  \mathbb{R}^n, 
\end{equation}
and 
\begin{equation} \label{omega-limit}
\lim_{\mu \rightarrow \infty} \omega_0(s, \mu) = 0, \ \text{uniformly in}\  \  s \in \overline{\mathbb{K}} .
\end{equation}
The focus of this paper is when the matrix $\tilde{\A}$ is singular, and hence so is the vector field coefficient $\A$. We decompose the  $\tilde{\A}$ into symmetric and skew-symmetric parts
\begin{equation} 
\tilde{\A}(x,t) = A(x,t)  +D(x,t)
\end{equation}
and assume that
\begin{equation} \label{ellipticity}
\left\{
\begin{split}
& A(x,t)^* = A(x,t), \quad D(x,t)^* = - D(x,t)^*, \quad \text{for a.e.} \quad (x,t) \in \Omega_T, \\
&  \norm{A}_{L^\infty(\Omega_T)} \leq \Lambda, \quad \|D\|_{L^\infty((0,T),\textup{BMO}(\R^n))} \leq M_0, \\
& \Lambda^{-1} |\xi|^2 \leq \wei{A(x,t) \xi, \xi}, \quad \forall \ \xi \in \R^n, \quad \text{for a.e.} \quad (x,t) \in \Omega_T,
\end{split}
\right.
\end{equation} 
where in the above $A^*, D^*$ denote the adjoint matrices of $A, D$ respectively, the $\textup{BMO}$-norm stands for the well-known John-Nirenberg semi-norm of functions of mean bounded oscillation.  Observe from \eqref{ellipticity} that $\tilde{\A} \not\in L^\infty(\Omega_T)$ due to the fact that its skew-symmetric part $D$ is just only assumed to be in $\textup{BMO}$.  For convenience in stating the results, we need a notation of the class of vector fields satisfying \eqref{Caratho-1}-\eqref{ellipticity}.
\begin{definition} Given an open bounded set $Q \subset \mathbb{R}^{n +1}$, an open interval $\mathbb{K} \subset \mathbb{R}$, and given numbers $\Lambda >0, M_0, M_1 >0$.  Let $\omega_0 : \K \times [0, \infty) \rightarrow [0, \infty)$ be continuous such that \eqref{omega-limit} holds and $\norm{\omega_0}_{\infty} \leq M_1$.  We denote $\mathbb{U}_{Q, \mathbb{K}}(\Lambda, M_0, M_1, \omega_0)$ be the set consisting all vector fields $\A: Q \times \mathbb{K} \times \mathbb{R}^n \rightarrow \mathbb{R}^n $ such that \eqref{Caratho-1}-\eqref{ellipticity} hold when replacing $\Omega_T$ by $Q$. Moreover, with a given $\A \in \mathbb{U}_{Q, \mathbb{K}}(\Lambda, M_0, M_1, \omega_0)$, the matrix $\tA$ defined in \eqref{asymp-Uh} is called asymptotical matrix of $\A$.
\end{definition} 

For any integrable function $f$ on a measurable $E \subset \mathbb{R}^n$, and for a positive Borel measure $\mu$, we denote  
\[
|E| = \int_{E} dx, \quad \mu(E) = \int_{E} \mu(x) dx, \quad \text{and} \quad \wei{f}_{E} = \fint_{E} f(y) dy = \frac{1}{|E|} \int_{E} f(y) dy.
\]
For each $\rho>0$ and $z_0 =(x_0, t_0) \in \R^{n+1}$, we denote the parabolic cylinder by
\[
Q_\rho(z_0) = B_\rho(x_0) \times \Gamma_\rho(t_0), \quad \text{where} \quad \Gamma_\rho(t_0) = (t_0 -\rho_0^2, t_0],
\]
and $B_\rho(x_0)$ is the ball in $\R^n$ centered at $x_0$ with radius $\rho$.  At this moment, we refer the readers to the definition of $A_p$ classes of Muchkenhoupt weights, and definitions of weak solutions in Section \ref{Preli-sec}. The following theorem is the main result of the paper.
\begin{theorem}  \label{W-1p-reg}  Let $\Lambda >0 , M_0 >0 , M_1 >0, M_2 \geq 1$ and $p > 2$. Then there are $\delta = \delta (\Lambda, M_0, M_1, M_2, p, n)$ sufficiently small and $\beta = \beta(\Lambda, M_0, n) >1$ such that  the following holds: Suppose that $\Omega$ is a $C^1$ domain, $\K \subset \R$ is an open interval and $\omega_0:\K\times [0, \infty)\rightarrow [0~,~\infty)$ is continuous satisfying \eqref{omega-limit} and $\norm{\omega_0}_{\infty} \leq M_1$.  Suppose also that $\A \in \mathbb{U}_{\Omega_T, \mathbb{K}} (\Lambda, M_0, M_1, \omega_0)$ with its asymptotical matrix $\tA$ satisfying \eqref{ellipticity} and 
\begin{equation}  \label{smallness-BMO-main}
[[\tA]]_{\textup{BMO}^\#(\Omega_T, \beta, R_0)} := \left( \sup_{0 < r < R_0} \sup_{z_0 = (x_0,t_0) \in \Omega_T} \frac{1}{|Q_r(z_0)|} \int_{Q_{r}(z_0) \cap \Omega_T} \Big| \tA(x,t) -\wei{\tA}_{B_r(x_0) \cap \Omega}(t) \Big|^\beta dxdt \right)^{1/\beta} \leq \delta, 
\end{equation}
for some $R_0>0$. Then, if $u \in L^2((0, T), W^{1,2}_0(\Omega))$ is a weak solution of  \eqref{main-eqn}, the estimate
\[
\begin{split}
\frac{1}{\omega(\Omega \times (\bar{t}, T))} \int_{\Omega \times(\bar{t}, T)}   |\nabla u|^p \omega(x,t) dx dt &  \leq C \left[ \frac{1}{\omega(\Omega \times(0, T) )} \int_{\Omega \times(0, T) } | \F|^p \omega(x, t) dx dt \right.\\
& \quad \quad + \left.  \left ( \fint_{\Omega \times(0, T)} |\nabla u|^2 dx dt\right)^{p/2}  +1 \right ]
\end{split}
\]
holds for $\bar{t} \in (0, T)$, for $\omega \in A_{p/2}$ with $[\omega]_{A_{p/2}} \leq M_2$, for some uniform constant $C$ depending only on $\Lambda, M_0, M_1, M_2, \Omega, R_0, \bar{t}, n, p$, and assuming $|\F| \in L^2(\Omega_T) \cap L^p(\Omega_T, \omega)$.
\end{theorem}
\noindent
A few comments on Theorem \ref{W-1p-reg} are now in ordered. Firstly, note that $\tA$ is not bounded, because we only assume that $D \in L^\infty((0,T), \textup{BMO})$. Therefore, even in linear and unweighted case, Theorem \ref{W-1p-reg} is already new. Note that since $D$ is skew-symmetric, it follows from \cite{CLMS}, see also that \cite{MV}, that for all $t \in (0,T)$ the following bilinear form is bounded: 
\[
\left|\int_{\Omega} \wei{D(x,t)\nabla u, \nabla v} dx \right| \leq C(n) \norm{D}_{L^\infty((0,T), \textup{BMO})} \norm{\nabla u}_{L^2(\Omega)} \norm{\nabla v}_{L^2(\Omega)}, \quad \forall u, v \in W^{1,2}_0(\Omega).
\]
This is one of the key ingredients for our paper. Observe also that the skew-symmetric part $D$ can be derived from the divergence-free drift term, see \cite{Sverak}. Therefore, the equations \eqref{main-eqn} can be rewritten as equations with singular divergence-free drifts. This class of equations are of great interest and has been investigated in many papers due to its relevance in many applications such as in fluid dynamics, and biology, see \cite{KS, LS, TP-note, Sverak, Q-Zhang} to cite a few. 

Secondly, we emphasize that  Theorem \ref{W-1p-reg} treats the case that coefficients $\A$ are unbounded and they could depend on $u$.  Moreover, we do not require the continuity of $\A$ in $u$. We also refer the readers to \cite{BW-no, BW1, BW2, CP, CMP, CFL, Ragusa, Fazio-1, Kim-Krylov, KZ1, Krylov, Softova, MP-1, MP, Meyers, Trud, Wang} for other papers in the same directions but only for linear equations with uniformly elliptic, bounded coefficients or for nonlinear equations in which $\A$ is independent on $u$. Indeed, it is well-known that the establishment of theory of Calder\'on-Zygmund estimates relies heavily on scaling invariant, see \cite{Wang} for the geometric intuition.  As a consequence, Calder\'on-Zygmund estimates are intrinsically invariant under the dilation $u \rightarrow u_\lambda := u(\lambda x, \lambda^2 t)/\lambda$ with $\lambda >0$,  see the estimates \eqref{CZ-est}, \eqref{W-1p-flat} and Remark \ref{remark-2} at the end of the paper. Therefore, the Calder\'on-Zygmund estimates are usually available only for the PDEs which are invariant under this dilation. For example, see \cite{BW-no, BW1, BW2, CP, CMP, CFL, Ragusa, Fazio-1, Kim-Krylov, KZ1, Krylov, Softova, MP-1, MP, Meyers, Trud, Wang} for which linear equations and  nonlinear equations where $\A$ is independent on $u$ are studied, and those equations are invariant under the dilations $u \rightarrow u_\lambda$. However, as $\A$ depends on $u$, the equation \eqref{main-eqn} will be changed under the dilations $u \rightarrow u_\lambda$ and this creates a serious issue. Only very recently are there  a few papers treating this case, see \cite{Bolegein, LTT, TN, NP}. In this paper, we overcome the inhomogeneity in \eqref{main-eqn} by adapting  the perturbation technique with double-scaling parameter method introduced in \cite{LTT}. See also \cite{TN, TP-degen, NP} for the implementation of the method. We essentially enlarge and consider the following class of equations with scaling parameter
\begin{equation} \label{lambda-class-eqn}
u_t - \text{div}[\A_\lambda(x,t, u, \nabla u)] = \text{div}[\F],
\end{equation}
where $\lambda >0$ is the scaling parameter, and 
\begin{equation} \label{A-lambda}
 \A_{\lambda}(x,t, s,\eta) = \A(x, t, \lambda s, \lambda \eta)/\lambda, \quad (x, t, s, \eta) \in \Omega_T \times \K_\lambda \times \mathbb{R}^n.
\end{equation}
Observe that \eqref{main-eqn} is just a special case of the \eqref{lambda-class-eqn}. We then prove a version of Theorem \ref{W-1p-reg} for \eqref{lambda-class-eqn} for all positive $\lambda$ and then obtain Theorem \ref{W-1p-reg} as a special case. In this perspective, the following observation regarding the scaling property of \eqref{main-eqn} is essential in the paper.
\begin{remark} \label{remark-1} Let  $\A \in \U_{\Omega_T, \K}(\Lambda, M_0, M_1, \omega_0)$ with its corresponding asymptotical matrix $\tA$. Let $\lambda >0$ and define $\K_\lambda  = \K/\lambda$. 
Then, for $\A_\lambda$ defined as in \eqref{A-lambda}, it is simple to check that 
\begin{equation} \label{lambda-asymp-Uh}
|\A_\lambda(x,t, s,\eta) - \tA(x,t) \eta| \leq \frac{1}{\lambda}\omega_0(\lambda s, \lambda \eta) (1 + |\lambda \eta|), \quad 
\forall \ (x,t, s,\eta) \in \Omega_T \times \K_\lambda \times \R^n.
\end{equation}
\end{remark} 

Lastly, we would like to emphasize that the setting that $\A(x,t, \cdot , \eta)$ is only defined on $ \K$, a  subset $\R$, is important in many applications, see \cite{LTT} for an example for which $\K = (0, K)$ with some positive $K$. Note also that when $\tA$ is bounded,
\[
[[\tA]]_{\textup{BMO}^\#(\Omega_T, \beta, r_0)} \sim [[\tA]]_{\textup{BMO}^\#(\Omega_T, 1, r_0)},\quad \forall \beta>1.
\]
Hence, the smallness condition \eqref{smallness-BMO-main} is reduced to the one required in \cite{BW-no, BW1, BW2, CP, CMP, CFL, Ragusa, Fazio-1, Kim-Krylov, KZ1, Krylov, Softova, MP-1, MP, Trud, Wang} for bounded coefficient $\tA$. Therefore, Theorem \ref{W-1p-reg} covers the results in these papers.  Moreover, it is known that this smallness condition is optimal as there is counterexample provided in \cite{Meyers} for uniformly elliptic, bounded coefficients. 

Besides global regularity estimates as in Theorem \ref{W-1p-reg}, we  are also interested in the theory of local regularity estimates. This is due to the fact that the local regularity estimates are sometimes important in applications as they only requires local information on the given data. Moreover, local regularity estimates are known to give the global ones.  Therefore, on one hand, it worths studying and stating local regularity estimates separately, and explicitly. On the other hand, note that local regularity theories such as Theorems \ref{interior-thm}-\ref{W-1p-reg-b} below in general could  not be derived from the global ones. Moreover, note that only global regularity estimates are available for the type of non-smooth domains studied in \cite{BW1, BW2, MP-1, MP}. Due to our interest, we therefore only consider $C^1$-domains in this paper, but non-smooth domains could be investigated in some projects in near future.

We next state our results for local Calder\'on-Zygmund type regularity estimates. When $z_0 = 0$, we omit it and write $Q_\rho = Q_\rho(0)$ for every $\rho >0$.  The following interior regularity estimate is our next result in the paper.
\begin{theorem} \label{interior-thm}  Let $\Lambda >0 , M_0 >0 , M_1 >0, M_2 \geq 1$ and $p > 2$. Then there are $\delta = \delta (\Lambda, M_0, M_1, M_2, p, n)$ sufficiently small and 
$\beta = \beta(\Lambda, M_0, n) >1$ such that  the following holds: Suppose that $\K \subset \R$ is an open interval and $\omega_0:\K\times [0, \infty)\rightarrow [0~,~\infty)$ is continuous satisfying \eqref{omega-limit} and $\norm{\omega_0}_{\infty} \leq M_1$.  For some $R>0$, let $\A \in \mathbb{U}_{Q_{2R}, \mathbb{K}} (\Lambda, M_0, M_1, \omega_0)$, with its asymptotical matrix $\tA$ satisfying \eqref{ellipticity} on $Q_{2R}$ and
\begin{equation*}
[[\tA]]_{\textup{BMO}^\#(Q_R, \beta, R)}  : = \left( \sup_{0 < \rho < R} \sup_{z_0 = (x_0, t_0)  \in Q_R} \fint_{Q_\rho(z_0)} |\tA(x,t) - \wei{\tA}_{B_\rho(x_0)}(t)|^\beta dxdt \right)^{1/\beta} \leq \delta.
\end{equation*}
Then, for every weak solution $u \in L^2(\Gamma_{2R}, W^{1,2}(B_{2R}))$ of 
\[
u_t - \textup{div}[\A(x,t, u, \nabla u)] = \textup{div}[\F], \quad \text{in} \quad Q_{2R},
\]
it holds that
\begin{equation} \label{CZ-est}
\begin{split}
& \left(\frac{1}{\omega(Q_{R})} \int_{Q_{R}} |\nabla u(x,t)|^p \omega(x,t) dx dt \right)^{1/p} \\
& \leq C\left[ \left(
\frac{1}{|Q_{2R}|}\int_{Q_{2R}}|\nabla u|^2 dxdt  \right)^{1/2} + \left(\frac{1}{\omega(Q_{2R})}\int_{B_{2R}} | \F |^p \omega(x, t) dx dt \right)^{1/p}  +1\right],
\end{split}
\end{equation}
 assuming that $|\F| \in L^2(Q_{2R}, \omega) \cap L^p(Q_{2R}, \omega)$, for $\omega \in A_{p/2}$ satisfying $[\omega]_{A_{p/2}} \leq M_2$, where $C$ is a constant dependent on $\Lambda, M_0, M_1, M_2, \omega_0, p, n$.
\end{theorem}
\noindent
Our last theorem is about regularity estimates on the flat domain, for which we define the upper balls in $\R^n$ to be
\[
B_R^+ = \{x = (x', x_n) \in B_R: x_n >0\}, \quad T_R = \{ x = (x', x_n) \in B_R : x_n =0\}.
\]
We also denote $Q_{R}^+ = B_R^+ \times \Gamma_R$. We then can state our result as below.
\begin{theorem} \label{W-1p-reg-b}  Let $\Lambda >0 , M_0 >0 , M_1 >0, M_2 \geq 1$ and $p > 2$. Then there are $\delta = \delta (\Lambda, M_0, M_1, M_2, p, n)$ sufficiently small and $\beta = \beta(\Lambda, M_0, n) >1$ such that  the following holds: Suppose that $R>0$, $\K \subset \R$ is an open interval and $\omega_0:\K\times [0, \infty)\rightarrow [0~,~\infty)$ is continuous satisfying \eqref{omega-limit} and $\norm{\omega_0}_{\infty} \leq M_1$.  Suppose also that $\A \in \mathbb{U}_{Q_{2R}^+, \mathbb{K}} (\Lambda, M_0, M_1, \omega_0)$, with its asymptotical matrix $\tA$ satisfying \eqref{ellipticity} on $Q_{2R}^+$ and 
\begin{equation*} 
[[\tA]]_{\textup{BMO}\#(Q_{R}^+, \beta, R)} : = \left( \sup_{0 < \rho < R} \sup_{z_0 = (x_0, t_0)  \in Q_R^+} \fint_{Q_\rho(z_0) \cap Q_{2R}^+}  |\tA(x,t) - \wei{\tA}_{B_\rho(x_0) \cap B_{2R}^+}(t)|^\beta dxdt \right)^{1/\beta} \leq \delta.
\end{equation*}
Then, for every weak solution $u \in L^2(\Gamma_{2R}, W^{1,2}(B_{2R}^+))$ of
\[
\left\{
\begin{array}{cccl}
u_t - \textup{div}[\A(x,t, u,\nabla u)] & = & \textup{div}[\F], & \quad \text{in} \quad Q_{2R}^+, \\
u & =& 0, & \quad \text{on} \quad T_{2R} \times \Gamma_{2R},
\end{array} \right.
\]
the following estimate holds
\[
\begin{split}
\left( \frac{1}{\omega(Q_R^+)}\int_{Q_{R}^+}   |\nabla u|^p \omega(x,t) dx dt \right)^{1/p}& \leq C \left[ \left( \frac{1}{\omega(Q_{2R}^+)} \int_{Q_{2R}^+} | \F|^p \omega(x, t) dx dt \right)^{1/p}\right. \\
& \quad \quad + \left.  \left( \fint_{Q_{2R}^+} |\nabla u|^2 dx dt \right)^{1/2} + 1 \right]
\end{split}
\]
for some uniform constant $C$ depending only on $\Lambda, M_0, M_1, M_2, \omega_0, n, p$, and for some $\omega \in A_{p/2}$ with $[\omega]_{A_{p/2}} \leq M_2$, assuming that $|\F| \in L^2(Q_{2R}^+) \cap L^p(Q_{2R}^+, \omega)$.
\end{theorem}

We now conclude the section by highlighting the layout of the paper.  Definitions of weak solutions, other definitions,  some analysis preliminary tools in measure theories and weighted norm inequalities are reviewed in the next section, Section \ref{Preli-sec}. The existence and uniqueness of weak solutions of linear equations with singular coefficients are also proved in this section, Theorem \ref{existence-theorem}.   Section \ref{interior-approx-section} consists interior intermediate step estimates and the proof of Theorem \ref{interior-thm}.  In fact, this section states and proves Theorem \ref{interior-rg-thm}, a more general version of Theorem \ref{interior-thm}. Section \ref{global-regularity-section} treats the estimates on the flat boundary and then states and proves Theorem \ref{main-theorem}, a slightly more general version of Theorem \ref{W-1p-reg-b}. Proof of Theorem \ref{W-1p-reg} will be provided in the last section, Section \ref{main-theorm-proof}, of the paper.
\section{Definitions and Preliminary results} \label{Preli-sec}
\subsection{Existence and uniqueness of weak solutions for singular parabolic equations}
For each open, bounded domain $\Omega \subset \R^n$ with Lipschitz boundary $\partial \Omega$. Let us denote $C_0^\infty(\Omega)$ the set of all smooth, compactly supported functions in $\Omega$, and $W^{1,2}_0(\Omega)$ is the closure of $C_0^\infty(\Omega)$ with the Dirichlet norm $\norm{u}_{W^{1,2}_0(\Omega)} = \norm{\nabla u}_{L^2(\Omega)}$. 
As usual $W^{-1,2}_0(\Omega)$ is the dual space of $W^{1,2}_0(\Omega)$. 
We also denote, 
\[
C_{0, p}^\infty(\Omega_T) = \left\{\varphi \in C_0^\infty(\overline{\Omega}_T): \varphi_{|\partial_p \Omega_T} =0 \right\}.
\] 
Also, let 
\[
\mathcal{E}_0(\Omega_T) = \text{closure of} \ C_{0, p}^\infty(\Omega_T) \ \text{in} \ L^2((0,T), W^{1,2}(\Omega)).
\]
\begin{definition} \label{weak-sol-def} Let $\Lambda >0, M_0$ be fixed and let $\tA$ such that \eqref{ellipticity} holds in $\Omega_T$. Also let  $\F  = (\F_k)_{k =1, 2,\cdot, n}$ be in $ L^2(\Omega_T)^{n}$,  and $g \in L^2((0,T), W^{1,2}(\Omega))$. 
\begin{itemize}
\item[\textup{(i)}] The function $u  \in L^2((0,T), W^{1,2}_{0}(\Omega)) \cap L^\infty((0,T), L^2(\Omega))$  is called a weak solution of the equations
\begin{equation} \label{u-eqn-def}
\left\{
\begin{array}{cccl}
u_t-  \textup{div}(\tA(x,t) \nabla u ) & = &  \textup{div}(\F), & \quad \text{in}\  \Omega_T, \\
u & =& g, & \quad \text{on} \quad \partial_p \Omega_T,
\end{array} \right.
\end{equation}
if $u_t \in L^2((0,T), W^{-1, 2}_{0}(\Omega), u -g  \in \mathcal{E}_0(\Omega_T)$, and
\begin{equation*} \label{weak-fornula}
\begin{split}
& \int_0^T \wei{u_t(\cdot, t), \varphi(\cdot, t)}_{W^{-1, 2}_{0}(\Omega), W^{1, 2}_{0}(\Omega)}  dt + \int_{\Omega_T} \wei{\tA (x,t)\nabla u (x,t), \nabla \varphi (x,t) } dxdt \\
& = -\int_{\Omega} \wei{\F(x,t), \nabla \varphi (x,t)} dx dt, \quad \forall \ \varphi \in C_{0,p}^\infty(\Omega_T).
\end{split}
\end{equation*}
\item[\textup{(iii)}] The function  $u  \in L^2((0,T), W^{1,2}(\Omega)) \cap L^\infty((0,T), L^2(\Omega))$ is called a weak solution of 
\begin{equation} \label{u-eqn-def-2}
\begin{array}{cccl}
u_t-  \textup{div}(\tA(x,t) \nabla u )  & = &  \textup{div}(\F), & \quad \text{in}\  \Omega_T
\end{array} 
\end{equation}
if $u_t \in L^2((0,T), W^{-1, 2}_{0}(\Omega)$, and  for all $\varphi \in C_{0,p}^\infty(\Omega_T)$,
\begin{equation*} \label{weak-fornula}
\int_0^T \wei{u_t(\cdot, t), \varphi(\cdot, t)}_{W^{-1, 2}_{0}(\Omega), W^{1, 2}_{0}(\Omega)} dt + \int_{\Omega_T} \wei{\tA (x,t)\nabla u (x, t), \nabla \varphi (x,t) } dx dt= -\int_{\Omega_T} \wei{\F(x, t) , \nabla \varphi (x,t)} dx dt.
\end{equation*}
\end{itemize}
\end{definition} 
\noindent
The following remark is important in this paper.
\begin{remark} \label{remark-MV} Observe that since $\tA = A + D$, with $A \in L^\infty(\Omega_T)$ and $D \in L^\infty((0,T), \textup{BMO})$ as in \eqref{ellipticity}, by \cite{CLMS, MV} we see that for  all $\phi \in W^{1,2}(\Omega)$ and $\varphi \in W^{1,2}_0(\Omega),
$
\[
\left| \int_\Omega \wei{\tA(\cdot, t) \nabla \phi, \nabla \varphi} dx\right| \leq C( \Lambda, M_0)\norm{\nabla \phi}_{L^2(\Omega)} \norm{\nabla \varphi}_{L^2(\Omega)}.
\]
Therefore, the term in the left hand-side of \eqref{weak-fornula} is well defined.
\end{remark}

The main result of this section is following existence, uniqueness theorem, which will be used frequently in our approximation estimates in Section \ref{interior-approx-section} and Section \ref{global-regularity-section}.
\begin{theorem}  \label{existence-theorem} Let $\Lambda, M_0, T$ be positive numbers and let $\Omega$ be a Lipschitz domain in $\R^n$. Assume that $\tA$ satisfies \eqref{ellipticity} in $\Omega_T$. Then for every $f \in L^2((0,T), W^{-1,2}_0(\Omega)$, there exists unique weak solution $u \in L^2((0,T), W^{1,2}_{0}(\Omega)) \times L^\infty((0,T), L^2(\Omega))$ of 
\[
\left\{
\begin{array}{cccl}
u_t-  \textup{div}(\tA(x,t) \nabla u ) & = &  f, & \quad \text{in}\  \Omega_T, \\
u & =& 0, & \quad \text{on} \quad \partial_p \Omega_T,
\end{array} \right.
\]
Moreover, 
\begin{equation} \label{energy-existence}
\sup_{t \in (0,T)}\norm{u}_{L^2(\Omega)} + \norm{\nabla u}_{L^2(\Omega_T)}  + \norm{u_t}_{L^2((0,T), W^{-1,2}_{0}(\Omega))} \leq C(\Lambda, M_0)  \norm{f}_{L^2((0,T), W^{-1,2}_0(\Omega))}.
\end{equation}
\end{theorem}
\begin{proof} Since our coefficients $\tA$ is singular, which has not studied elsewhere,  a proof for Theorem \ref{existence-theorem} is needed. We use Galerkin's approximation method as in \cite[p. 353-358]{Evans} with some modification. We only outline some main steps. Let us define the following bilinear map
\[
\begin{split}
& B : W^{1,2}_{0}(\Omega) \times W^{1,2}_{0}(\Omega) \rightarrow \R, \quad \text{with} \\
& B(u,v; t) =  \int_{\Omega} \wei{\tA(x, t) \nabla u(x), \nabla v(x)} dx, \quad \forall u, \ v  \in W^{1,2}_{0}(\Omega).
\end{split}
\]
It follows from Remark \ref{remark-MV} that
\begin{equation} \label{bounded-B}
|B(u, v;t)| \leq C(\Lambda, M_0)\norm{\nabla u}_{L^2(\Omega)} \norm{\nabla v}_{L^2(\Omega)}, \quad \forall \ u, v \in W^{1,2}_{0}(\Omega).
\end{equation}
Moreover, due to the assumption that $\tA = A +D$, and $D$ is skew-symmetric, and due to the assumption on ellipticity of $A$ in \eqref{ellipticity}, we observe that
\begin{equation} \label{coercive-B}
B(u,u, t) = \int_\Omega \wei{A(x,t) \nabla u, \nabla u} dx  \geq \Lambda \norm{\nabla u}^2_{L^2(\Omega)}.
\end{equation}
Therefore, $B$ is bilinear, bounded, and coercive.  Now, let $\{w_k\}_{k=1}^\infty$ be in $C^\infty_{0}(\Omega)$. Moreover, $\{w_k\}_{k=1}^\infty$  is an orthogonal basis of $W^{1,2}_{0}(\Omega)$, and an orthonormal basic of $L^2(\Omega)$. For each $m \in \N$, we look for the approximation solution $\bu_m : [0, T] \rightarrow W^{1,2}_{0}(\Omega)$ of the form
\[
\bu_m(t) = \sum_{k=1}^m d_m^k(t) w_k,
\]
where, $d_m^k$ satisfies the equation
\begin{equation} \label{ODE-d}
\left\{
\begin{array}{cccl}
\partial_t d_m^k(t) + \displaystyle{\sum_{k=1}^m e_{kl}(t) d_m^l(t)} & = & f^k(t), & \quad t \in (0,T), \\
d_m^k(0) & = & 0, & \quad k = 1,2,\cdots, m,
\end{array} \right.
\end{equation}
with
\[
e_{kl}(t) = B(w_l, w_k;t), \quad f_k(t) = \wei{f(\cdot, t), w_k}_{W^{-1,2}_0(\Omega), W^{1,2}_0(\Omega)}. 
\]
From \eqref{bounded-B}, we see that
\[
\norm{e_{kl}}_{L^\infty(0,T)} \leq C(\Lambda, M_0), \quad \norm{f_k}_{L^2((0,T)} \leq \norm{f}_{L^2((0,T), W^{-1,2}_0(\Omega))} \norm{w_k}_{W^{1,2}_0(\Omega)}, \quad \forall \ k, l = 1, \cdots, m.
\]
Hence, the existence of solutions for the system \eqref{ODE-d} follows by standard ODE theory. Moreover, using \eqref{bounded-B}-\eqref{coercive-B}, we can follow the energy estimates as in \cite[p. 353-358]{Evans} to obtain
\begin{equation*} 
\sup_{t\in (0,T)} \norm{\bu_m}_{L^2(\Omega)} + \norm{\bu_m}_{L^2((0,T), W^{1,2}_{0}(\Omega))} + \norm{\partial_t \bu_m}_{L^2((0,T), W^{-1,2}_{0}(\Omega))} \leq C(\Lambda, M_0) \leq \norm{f}_{L^2((0,T), W^{-1,2}_0(\Omega))}.
\end{equation*}
From this estimate, and as in \cite[p. 353-358]{Evans}, we can pass through limit as $m \rightarrow \infty$ to obtain the existence of $u \in L^2((0,T), W^{1,2}_{0}(\Omega)) \cap L^\infty((0,T), L^2(\Omega))$ satisfying the estimate
\begin{equation} \label{u-existence-energy}
\sup_{t\in (0,T)} \norm{u}_{L^2(\Omega)} + \norm{u}_{L^2((0,T), W^{1,2}_{0}(\Omega))} + \norm{u_t}_{L^2((0,T), W^{-1,2}_{0}(\Omega))} \leq C(\Lambda, M_0) \norm{f}_{L^2((0,T), W^{-1,2}_0(\Omega))}.
\end{equation}
Moreover, $u \in \mathcal{E}_0(\Omega_T)$ and for every $ \varphi \in C^{\infty}_{0,p}(\Omega_T)$
\begin{equation*}
\int_0^T \wei{u_t(\cdot, t), \varphi(\cdot, t)}_{W^{-1,2}_{0}(\Omega), W^{1,2}_{0}(\Omega)} dt +  \int_0^T B(u(\cdot, t), \varphi(\cdot, t); t) dt + \int_{\Omega_T} \wei{\F(x, t), \nabla \varphi(x,t)} dxdt =0.
\end{equation*}
The uniqueness of solutions also follows from \eqref{energy-existence} and the linearity of our considered equations.  This completes the proof of the theorem.
\end{proof}
\subsection{Gehring's type regularity estimates for singular homogeneous equations}
This section states two self-improved regularity estimates for weak solutions of our class of singular parabolic equations, assuming that the skew-symmetric part $D \in L^\infty_t(\textup{BMO})$. These type estimates are sometimes referred  as Meyer's type estimates, see \cite{Gehring, Meyers}. The results in this section are new, and of independent interests. They also improve the classical results for which only the case $D =0$ is studied (i.e. $\A = A$ is symmetric).
The first main result of the section is the following result.
\begin{lemma}[Gehring's type regularity]  \label{Gehring-1} Let $\Lambda, M_0$ be positive numbers. Then, there exists $\gamma = \gamma(\Lambda, M_0, n) >0$ such that the following statement holds. Assume that \eqref{ellipticity} holds for a given matrix $\tA$ in $Q_{7/4}$. If $v \in L^2(\Gamma_{7/4}, W^{1,2}(B_{7/4}))$ is a weak solution of 
\[
v_t - \textup{div}[\tA(x,t) \nabla v(x,t)] =0, \quad \text{in} \quad Q_{7/4},
\]
then
\[
\left(\fint_{Q_{3/2}} |\nabla v|^{2+\gamma} dxdt \right) ^{\frac{1}{2+\gamma}} \leq C(\Lambda, M_0, n) \left( \fint_{Q_{7/4}} |\nabla v|^2 dxdt \right)^{1/2}.
\]
\end{lemma}
\begin{proof}  We skip the proof because it is similar to that of \cite[Lemma 2.6]{TP-NS}.
\end{proof}
\noindent
Now, for $\rho >0$, we denote 
\[
B_\rho^+ = \Big \{x = (x', x_n) \in B_\rho: x_n > 0 \Big \}, \quad Q_\rho^+ = B_\rho^+ \times \Gamma_\rho, \quad \Gamma_\rho = (-\rho^2, 0]. 
\]
 Moreover, the flat boundary part of $B_\rho^+$ is denoted by $T_\rho$, i.e.
 \[
 T_\rho = \Big\{x = (x', 0) \in B_\rho \Big\}.
 \]
The following version of Gehring's type estimate on the flat boundary domain is also needed in the paper.
\begin{lemma} \label{Gehring-2} Let $\Lambda>0, M_0 >0$. Then, there exists $\gamma = \gamma(n, \Lambda, M_0) >0$ such that the following statement holds: Assume that $\tA$ satisfy the assumption \eqref{ellipticity} on $Q_{7/4}^+$. If $v  \in L^2(\Gamma_{7/4}, W^{1,2}_\sigma(B_{7/4}^+))$ is a weak solutions of
\begin{equation} \label{V-Gehring-1-eqn}
\left\{
\begin{array}{cccl}
v_t - \textup{div} [\tA(x,t) \nabla v]  &= & 0, & \quad \text{in} \quad Q_{7/4}^+, \\
v & = & 0, & \quad \text{on} \quad T_{7/4}\times \Gamma_{7/4},
\end{array} \right.
\end{equation}
then, 
\[
\left(\fint_{Q_{3/2}^+} |\nabla v|^{2+\delta} \right)^{1/(2+\delta)} \leq C(n, \Lambda, M_0) \left(\fint_{Q_{7/4}^+} |\nabla v|^2 dx \right)^{1/2}.
\]
\end{lemma} 
\begin{proof} The proof of this lemma is almost identical to that of Lemma \ref{Gehring-1}. We therefore skip it, see also \cite[Lemma 2.7]{TP-NS}.
\end{proof}
\subsection{Lipschitz regularity estimates for homogeneous constant coefficient equations} 

In this section, let $A_0 : \Gamma_{3/2} \rightarrow \R^{n\times n}$ be a measurable $n\times n$ symmetric matrix, bounded, satisfying 
\begin{equation} \label{A-0-ellip}
\Lambda |\xi|^2 \leq \wei{A_0(t) \xi, \xi}, \quad \norm{A_0}_{L^\infty(\Gamma_{3/2})} \leq \Lambda^{-1}, \quad \text{for a.e.} \ t \in \Gamma_{3/2}, \ \forall \ \xi \in \R^n.
\end{equation}
The following regularity lemma is a well-known result for linear parabolic equations.
\begin{lemma} \label{L-infty-reg-inter} Assume that  the $n\times n$ symmetric matrix $A_0$ satisfies \eqref{A-0-ellip}. Then, if $w \in L^2(\Gamma_{3/2}, W^{1,2}(B_{3/2}))$ is a weak solution of 
\[
\begin{array}{cccl}
w_t - \textup{div}[A_0(t) \nabla w]  & =& 0, & \quad \text{in} \quad Q_{3/2},
\end{array}
\]
it holds that
\[
\norm{\nabla w}_{L^\infty(Q_1)}   \leq C(\Lambda, n) \left( \fint_{Q_{3/2}} |\nabla w|^2 dx dt \right)^{1/2}.
\]
\end{lemma}
\noindent
Similarly, the following Lipschitz regularity estimates for weak solutions on the flat domain is also well-known.
\begin{lemma} \label{L-infty-reg-bdr} Assume that  the $n\times n$ symmetric matrix $A_0$ satisfies \eqref{A-0-ellip}. Then, if $w \in L^2(\Gamma_{3/2}, W^{1,2}(B_{3/2}^+))$ is a weak solution of 
\[
\left\{
\begin{array}{cccl}
w_t - \textup{div}[A_0(t) \nabla w]  & =& 0, & \quad \text{in} \quad Q_{3/2}^+, \\
w & =& 0, & \quad \text{on} \quad T_{3/2} \times \Gamma_{3/2},
\end{array}
\right.
\]
it holds that
\[
\norm{\nabla w}_{L^\infty(Q_1^+)}   \leq C(\Lambda, n)  \left( \fint_{Q_{3/2}^+} |\nabla w|^2 dx dt \right)^{1/2}. 
\]
\end{lemma}

\subsection{Munckenhoupt weights and weighted inequalities} \label{Muckenhoupt-w}

For each $1 \leq q < \infty$, a non-negative, locally integrable function $\mu :\R^{n+1} \rightarrow [0, \infty)$ is said to be in the class of parabolic $A_q$ Muckenhoupt weights if 
\[
\begin{split}
[\mu]_{A_q} & :=  \sup_{r >0, z \in \R^{n+1}} \left(\fint_{Q_r(z)} \mu (x,t) dx dt  \right) \left(\fint_{Q_r(z)} \mu (x,t)^{\frac{1}{1-q}} dx dt \right)^{q-1} < \infty, \quad \textup{if} \quad q > 1, \\
[\mu]_{A_1} &: =  \sup_{ r>0 , z \in \R^{n+1}} \left(\fint_{Q_r(z)} \mu (x, t) dx dt \right)  \norm{\mu^{-1}}_{L^\infty(Q_{r}(z))}  < \infty \quad \textup{if} \quad q  =1.
\end{split}
\]
It is well known that the class of $A_p$-weights satisfies the reverse H\"{o}lder's inequality and the doubling properties, see for example \cite{Coif-Feffer, Fabes-Riviere, Stein}. In particular, a measure with an $A_p$-weight density is, in some sense, comparable with the Lebesgue measure. 
\begin{lemma}[\cite{Coif-Feffer}] \label{doubling} For $1 < q < \infty$, the following statements hold true
\begin{itemize}
\item[\textup{(i)}] If $\mu \in A_{q}$,  then for every parabolic cube $Q \subset \R^{n+1}$ and every measurable set $E\subset Q$, 
$
\mu(Q) \leq [\mu]_{A_{q}} \left(\frac{|Q|}{|E|}\right)^{q} \mu(E).
$
\item[\textup{(ii)}] If $\mu \in A_q$, then there is $C = C([\mu]_{A_q}, n)$ and $\beta = \beta([\mu]_{A_q}, n) >0$ such that
$
\mu(E) \leq C \left(\frac{|E|}{|Q|} \right)^\beta \mu(Q),
$
 for every parabolic cube $Q \subset \R^{n+1}$ and every measurable set $E\subset Q$.
\end{itemize}
\end{lemma} \noindent
Let us also recall the definition of the parabolic Hardy-Littlewood maximal operators which will be needed in the paper
\begin{definition} \label{Maximal-Operator} The parabolic Hardy-Littlewood maximal function of a locally integrable function  $f$ on $\R^{n+1}$ is defined by
\[
 (\M f)(x,t) = \sup_{\rho>0} \fint_{Q_\rho(x,t)}{|f(y,s)|\, dy ds}.
\]
If $f$ is defined in a domain $U\subset \R^n\times \R$, then we   denote
\[
 \M_{U} f = \M(\chi_U f).
\]
\end{definition}
\noindent
The following boundedness of the Hardy-Littlewood maximal operator is due to Muckenhout \cite{Muckenhoupt}. For the proof of this lemma, one can find it in \cite{Fabes-Riviere, Stein}.
\begin{lemma} \label{Hardy-Max-p-p} Assume that $\mu \in A_q$ for some $1 < q < \infty$. Then, the followings hold.
\begin{itemize}
\item[(i)] Strong $(q,q)$: There exists a constant $C = C([\mu]_{A_q}, n,q)$ such that  
\[ \|\mathcal{M} \|_{L^{q}(\R^{n+1},\mu) \to L^{q}(\R^{n+1}, \mu)} \leq C. \] 
\item[(ii)] Weak $(1,1)$: 
There exists a constant $C=C(n)$ such that for any $\lambda >0$, we have 
\[
\big| \big \{ (x,t) \in \mathbb{R}^{n+1}: \mathcal{M} (f) > \lambda \big\}\big| \leq \frac{C}{\lambda} \int_{\mathbb{R}^{n+1}}|f(x,t)| dxdt.
\]  
\end{itemize}
\end{lemma}
We also collect some useful measure theory results needed in the paper. Our first lemma is the standard result in in measure theory.
\begin{lemma} \label{measuretheory-lp}
Assume that $g\geq 0$ is a measurable function in a bounded subset $U\subset \mathbb{R}^{n+1}$. Let $\theta>0$ and $N>1$ be given constants. If $\mu$ is a weight in  $L^{1}_{loc}(\mathbb{R}^{n+1})$, then for any $1\leq p < \infty$ 
\[
g\in L^{p}(U,\mu) \Leftrightarrow S:= \sum_{j\geq 1} N^{pj}\mu(\{x\in U: g(x)>\theta N^{j}\}) < \infty. 
\]
Moreover, there exists a constant $C>0$ such that 
\[
C^{-1} S \leq \|g\|^{p}_{L^{p}(U,\mu)} \leq C (\mu(U) + S), 
\]
where $C$ depends only on $\theta, N$ and $p$. 
\end{lemma} \noindent
The following lemma is commonly used, and it is a consequence of the Vitali's covering lemma. The proof of this lemma can be found in \cite[Lemma 3.8]{MP-1}, which is in turn an extension of  \cite[Theorem 3]{Wang}. 
\begin{lemma}\label{Vitali} Let $\mu$ be an $A_{q}$ weight for some $q \in (1,\infty)$ with $[\mu]_{A_q} \leq M$ for some $M \geq 1$. Assume that $E \subset K \subset Q_1$ are measurable sets for which there exists $\epsilon, \rho_{0}\in (0, 1/4)$ such that 
\begin{itemize}
\item[(i)]  $\mu(E) < \epsilon \mu(Q_{1}(z)) $ for all $z\in \overline{Q}_{1}$, and 
\item[(ii)] for all $z \in Q_{1}$ and $\rho \in (0, \rho_{0}]$, if $\mu (E \cap Q_{\rho}(z)) \geq \epsilon \mu(Q_{\rho}(z))$, 
then $Q_{\rho}(z)\cap Q_{1} \subset K. $
\end{itemize}  
Then with $\e_1 = \e  (20)^{(n+2)q}M^2 $, the following estimate holds 
\[
\mu(E) \leq \epsilon_1 \,\mu(K). 
\] 
The same conclusion also holds if we replace $Q_1$ by $Q_1^+$.
\end{lemma}

\section{Interior weighted $W^{1,p}$- regularity theory} \label{interior-approx-section}
In this section, let $R>0$, and $\K \subset \R$ be an open interval. For $\A \in \U_{Q_{2R}, \K}(\Lambda, M_0, M_1, \omega_0)$, let $\A_\lambda$ be defined as in \eqref{A-lambda} with some $\lambda >0$. We focus on the following equation with scaling parameter $\lambda$:
\begin{equation} \label{eqn-in}
u_t - \text{div}[\A_\lambda(x, t, u, \nabla u)] = \text{div}[\F] , \quad \text{in} \quad Q_{2R}.
\end{equation}
Observe that the equation  \eqref{main-eqn} can be reduced from in \eqref{eqn-in} by taking $\lambda =1$. As we already discussed, enlarging the class of equations in \eqref{main-eqn} to the class of equations in \eqref{eqn-in} with the parameter $\lambda>0$ is essential in our approach due to the homogeneity in the Calder\'on-Zygmund type regularity estimates.

Let $\tA$ be the asymptotical matrix of $\A$, we recall that $\tA = A +D$, where $A$ is an $n\times n $ symmetric matrix, and $D$ is a $n\times n$ skew-symmetric matrix. Moreover, $A \in L^\infty(Q_{2R})$,  $D \in L^\infty(\Gamma_{2R}, \textup{BMO})$, and
\begin{equation} \label{ellip-interior}
\Lambda^{-1} |\eta|^2  \leq \wei{A(x,t) \eta, \eta}, \quad \norm{A}_{L^\infty(Q_{2R})} \leq \Lambda, \quad  
\norm{D}_{L^\infty(\Gamma_R, \textup{BMO})} \leq M_0, \quad \text{for a.e.} \ (x,t) \in Q_{2R}, \quad \forall \ \eta \in \R^n.
\end{equation}
Now, let  $\gamma = \gamma(\Lambda, M_0, n)>0$ be defined as in  Lemma \ref{Gehring-1}, and let 
\begin{equation} \label{beta.def}
\beta = \frac{2(2+\gamma)}{\gamma} > 1.
\end{equation}
For some $r_0 \in (0,R)$,  we denote
\[
[[\tA]]_{\textup{BMO}^\#(Q_R, \beta, r_0)} = \left( \sup_{z_0 = (x_0, t_0) \in Q_R} \sup_{0 < \rho < r_0} \fint_{Q_{\rho}(z_0)} |\tA(x,t) -\wei{\tA}_{B_\rho(x_0)}(t)|^\beta dxdt \right)^{1/\beta},
\]
where $\wei{\tA}_{B_\rho(x_0)}(t)$ denotes the mean of $\tA$ on $B_\rho(x_0)$. This section proves the following theorem, which in turns also produces Theorem \ref{interior-thm} when taking $\lambda =1$.
\begin{theorem} \label{interior-rg-thm} Let $\Lambda >0 , M_0 >0 , M_1 >0, M_2 \geq 1$ and $p > 2$. Then there are $\delta = \delta (\Lambda, M_0, M_1, M_2, p, n)$ sufficiently small and 
$\beta = \beta(\Lambda, M_0, n) >1$ such that  the following holds: Suppose that $\K \subset \R$ is an open interval and $\omega_0:\K\times [0, \infty)\rightarrow [0~,~\infty)$ is continuous satisfying \eqref{omega-limit} and $\norm{\omega_0}_{\infty} \leq M_1$.  For some $R>0$, let $\A \in \mathbb{U}_{Q_{2R}, \mathbb{K}} (\Lambda, M_0, M_1, \omega_0)$, with its asymptotical matrix $\tA$ satisfying \eqref{ellip-interior} and
\begin{equation}\label{BMO-B-R}
[[\tA]]_{\textup{BMO}^\#(Q_R, \beta, R)}  \leq \delta.
\end{equation}
If $u \in L^2(\Gamma_{2R}, W^{1,2}(B_{2R}))$ 
is a weak solution of \eqref{eqn-in} with some $\lambda >0$, then
\begin{equation} \label{CZ-est}
\begin{split}
 \left(\frac{1}{\omega(Q_{R})} \int_{Q_{R}} |\nabla u(x,t)|^p \omega(x,t) dx dt \right)^{1/p} &  \leq C\left[\left(\frac{1}{\omega(Q_{2R})}\int_{B_{2R}} | \F |^p \omega(x, t) dx dt \right)^{1/p} \right.\\
 & \quad \quad \left. + \max\left\{ \left( \fint_{Q_{2R}}|\nabla u|^2 dxdt  \right)^{1/2}, \lambda^{-1} \right\} \right],
\end{split}
\end{equation}
holds, assuming that $|\F| \in L^2(Q_{2R}, \omega) \cap L^p(Q_{2R}, \omega)$, where $C = C(\Lambda, M_0, M_1, \omega_0, p, n)$, and $\omega \in A_{p/2}$ such that $[\omega]_{A_{p/2}} \leq M_2$.
\end{theorem}

\subsection{Interior approximation estimates} The following proposition is the main result of the subsection.
\begin{proposition} \label{interio-approx-lemma} Let $\Lambda >0, M_0 >0, M_1 >0$,  and let $\omega_0$ be a continuous function satisfying \eqref{omega-limit} and $\norm{\omega_0}_{\infty} \leq M_1$. Then, for every small number $\epsilon >0$, there exist $\delta = \delta(\epsilon, \Lambda, M_0, M_1, n) >0$ sufficiently small and sufficiently large $\lambda_0 =  \lambda_0(\epsilon, \Lambda, M_0, \omega_0, n) \geq 1$  such that the following holds:  Assume that $\A \in \U_{Q_{2R}, \K}( \Lambda, M_0, M_1, \omega_0)$, with its  asymptotical matrix $\tA$ satisfying \eqref{ellip-interior}, and for some $R>0$ and some open interval $\K \subset \mathbb{R}$. Assume also that
\[
[[\tA]]_{\textup{BMO}^\#(Q_R, \beta, r_0)} \leq \delta, \quad \text{for some} \quad r_0 \in (0,R).
\]
Then, for some $r \in (0,r_0/2)$ some $z_0 = (x_0, t_0) \in \overline{Q}_R$ and for $\lambda \geq \lambda_0$, if  
\[
\fint_{Q_{2r}(z_0)} | \F(x,t)|^2 dx dt  \leq \delta,
\]
and  if  $u$ is a weak solution of \eqref{eqn-in} satisfying
\[
\fint_{Q_{2r}(z_0)} |\nabla u|^2 dx dt \leq 1,
\]
then there is $v \in L^2(Q_{3r/2}(z_0))$, and constant $C = C(\Lambda, n)$ such that
\[
\fint_{Q_{3r/2}(z_0)} |\nabla u - \nabla v|^2 dx dt \leq \epsilon,
\quad \text{and} \quad \norm{\nabla v}_{L^\infty(Q_{r}(z_0))} \leq C.
\]
\end{proposition}
\begin{proof}  Without loss of generality, we assume assume $z_0 =0$ and $R =1$. We split the procedure for the proof into two steps of approximations.\\
{\bf Step I:} We  
approximate in $Q_{7r/4}$ the solution $u$ of \eqref{eqn-in} by the solution $w$ of the following equation
\begin{equation} \label{w-interio}
\left\{
\begin{array}{cccl}
w_t - \text{div}[\tilde{\A}(x,t) \nabla w] & =& 0, &\quad \text{in}\quad Q_{\frac{7r}{4}}, \\
w & =& u, & \quad \text{on} \quad \partial_p Q_{\frac{7r}{4}}.
\end{array}\right.
\end{equation}
Observe that for a given given $u \in L^2(\Gamma_{2r}, W^{1,2}(B_{2r})$ weak solution of \eqref{eqn-in}, $w$ is a weak solution of \eqref{w-interio} if and only if $\tilde{w} = w-u$ is a weak solution of
\begin{equation} \label{tilde-w-interior}
\left\{
\begin{array}{cccl}
\tilde{w}_t - \text{div}[\tilde{\A}(x,t) \nabla \tilde{w}] & =& \text{div}[G] &\quad \text{in}\quad Q_{\frac{7r}{4}}, \\
\tilde{w} & =& 0, & \quad \text{on} \quad \partial_p Q_{\frac{7r}{4}},
\end{array}\right.
\end{equation}
where $G = \tA(x,t) \nabla u - \A_\lambda(x,t, u, \nabla u) - \F$. It follows from \eqref{lambda-asymp-Uh} in Remark \ref{remark-1} that
\begin{equation} \label{A-differ}
|\tA(x,t) \nabla u - \A_\lambda(x,t, u, \nabla u)|  \leq \frac{1}{\lambda}  \omega_0( \lambda u, \lambda |\nabla u|) (1+ \lambda |\nabla u|).
\end{equation}
On the other hand, from \eqref{omega-limit} we see that for $\delta >0$, we can find a large number $K_\delta >0$ depending only on $\omega_0$ and $\delta$  such that $\omega_0(z, s) \leq \delta$ for all $s \geq K_\delta$ and for all $z \in \overline{\mathbb{K}} $. This and $\norm{\omega_0}_{\infty} \leq M_1$ in turn imply that
\begin{equation} \label{K-delta}
\omega_0(s, \mu) (1+\mu) \leq  \delta (1+\mu) + M_1 K_\delta, \quad \forall \mu \geq 0, \quad \forall s \in \K.
\end{equation}
Then, from this and \eqref{A-differ}, it follows that 
\begin{equation} \label{tA-A-lambda-diff}
|\tA(x,t) \nabla u - \A_\lambda(x,t, u, \nabla u)|  \leq \delta |\nabla u| + \lambda^{-1}(M_1 K_\delta + \delta).
\end{equation}
In particularly, by its definition, $G \in L^2(Q_{7r/4})$. From this, the existence of weak solution $\tilde{w}$ of \eqref{tilde-w-interior} follows from Theorem \ref{existence-theorem}. This implies the existence of $w$. Moreover, from Theorem \ref{existence-theorem}, and \eqref{tA-A-lambda-diff}, we also have
\begin{equation} \label{u-w-step-1}
\fint_{B_{\frac{7r}{4}}} |\nabla (u-w)|^2 dx dt  
\leq  C \left[ \delta^2 \fint_{B_{\frac{7r}{4}}} |\nabla u|^2 d x dt + \fint_{B_{\frac{7r}{4}}} | \F|^2dx dt   + \frac{(M_1 K_\delta +\delta)^2}{\lambda^2} \right],
\end{equation}
for some $C = C(\Lambda, M_0, n)$. On the other hand, it follows from \eqref{A-differ} that if $\lambda \geq 1$, then
\[
|G| \leq M_1(|\nabla u| +1) + |\F|.
\]
Therefore, we can also apply Theorem \ref{existence-theorem} again to yield
\[
\fint_{B_{\frac{7r}{4}}} |\nabla (u-w)|^2 dx dt  
\leq  C \left[ M_1^2 \fint_{B_{\frac{7r}{4}}} |\nabla u|^2 d x dt + \fint_{B_{\frac{7r}{4}}} | \F|^2dx dt   + M_1^2\right], \quad \text{if} \quad \lambda \geq 1.
\]
This estimate, and the assumptions in Proposition \ref{interio-approx-lemma} infer that
\begin{equation} \label{energy-w-step-1}
\begin{split}
\fint_{B_{\frac{7r}{4}}} |\nabla w|^2 dxdt  & \leq C(\Lambda, M_0, M_1, n), \quad \text{if} \quad \lambda \geq 1.
\end{split}
\end{equation}
Now, recall our choice of $\beta$ and $\gamma$ as in \eqref{beta.def}. From Lemma \ref{Gehring-1} and \eqref{energy-w-step-1}, we see that 
if $\lambda \geq 1$, then
\begin{equation} \label{Gehring-1-app}
\left[ \fint_{Q_{3r/2}(z_0)} |\nabla w|^{2+\gamma} dxdt \right]^{1/(2+\gamma)} \leq C(\Lambda, M_0, n) \left[\fint_{Q_{7r/4}(z_0)} |\nabla w|^2 dx dt\right ]^{1/2} \leq C(\Lambda, M_0, M_1, n).
\end{equation}
\noindent
{\bf Step 2:} In our second step, we approximate $w$ by $v$ which is a weak solution of the following equation
\begin{equation} \label{v-interior}
\left\{
\begin{array}{cccl}
v_t - \textup{div}[\wei{A}_{B_{r}}(t) \nabla v]  & = & 0, &  \quad \text{in} \quad Q_{3r/2}, \\
v & = & w, & \quad \text{on} \quad \partial_p Q_{3r/2}
\end{array} \right.
\end{equation}
Observe that since $\tA = A + D$, and  $\wei{D}_{B_r}(t)$ is skew-symmetric depending only on $t$, $v$ is a weak solution of \eqref{v-interior} if and only if it is 
a weak solution of 
\begin{equation} \label{v-eqn-inter-2}
\left\{
\begin{array}{cccl}
v_t - \textup{div}[\wei{\tA}_{B_{r}}(t) \nabla v]  & = & 0, &  \quad \text{in} \quad Q_{3r/2}, \\
v & = & w, & \quad \text{on} \quad \partial_p Q_{3r/2}
\end{array} \right.
\end{equation}
Let $\tilde{v} = v -w$. Then, $v$ is a weak solution of \eqref{v-eqn-inter-2} if and only if $\tilde{v}$ is a weak solution of
\begin{equation} \label{v-eqn-inter-3}
\left\{
\begin{array}{cccl}
\partial_t \tilde{v}_t - \textup{div}[\wei{\tA}_{B_r}(t) \nabla \tilde{v}]  & = &  \textup{div}[(\tA - \wei{\tA}_{B_r}(t)) \nabla w], &  \quad \text{in} \quad Q_{3r/2}, \\
\tilde{v} & = & 0, & \quad \text{on} \quad \partial_p Q_{3r/2}
\end{array} \right.
\end{equation}
We claim that $(\tA -\wei{\tA}_{B_{r}}(t)) \nabla w$ is in $ L^2(Q_{3r/2})$. Indeed, it follows from \eqref{Gehring-1-app}, our choice of $\beta$ as in \eqref{beta.def} and H\"{o}lder's inequality that
\begin{equation} \label{A-nabla-V-inter}
\begin{split}
 & \fint_{Q_{3r/2}} |(\tA -\wei{\tA}_{B_{r}}(t)) \nabla w|^2 dx dt\\
  & \leq \left(\fint_{Q_{3r/2}} |\nabla w|^{2+\gamma} dx dt \right)^{2/(2+\gamma)} \left(\fint_{Q_{r}}|\tA(x,t)-\wei{\tA}_{B_{r}}(t)|^{\frac{2(2+\gamma)}{\gamma}}  dxdt \right)^{\frac{\gamma}{(2+\gamma)}} \\
& \leq C(\Lambda, M_0, M_1, n) [[\tA]]_{\textup{BMO}^\#(Q_1, \beta, r_0)}^2.
\end{split}
\end{equation}
From the estimate \eqref{A-nabla-V-inter}, we can use Theorem \ref{existence-theorem} to obtain the existence of a weak solution $\tilde{v}$ of \eqref{v-eqn-inter-3}. Then, the existence of $v$ follows. Moreover, Theorem \ref{existence-theorem}  and \eqref{A-nabla-V-inter} also yield
\begin{equation*} \label{w-v-interior}
\fint_{Q_{3r/2}} |\nabla v -\nabla w|^2 dxdt \leq C(\Lambda, M_0, M_1, n) [[\tA]]_{\textup{BMO}^\#(Q_1, \beta, r_0)}^2, \quad \text{when} \quad \lambda \geq 1.
\end{equation*}
Now, we combine this estimate with \eqref{u-w-step-1} to infer that there is $C_* = C_*(\Lambda, M_0, M_1, n) >0$ such that if $\lambda \geq 1$, then
\begin{equation} \label{u-v-interior-final}
\fint_{Q_{3r/2}} |\nabla  u - \nabla  v|^2 dxdt \leq C_*(\Lambda, M_0, M_1, n) \Big[ \delta  +  \frac{(M_1 K_\delta +\delta)^2}{\lambda^2} \Big].
\end{equation}
Now, choose $\delta = \delta(\epsilon, \Lambda, M_0, M_1, n)>0$ sufficiently small such that 
\[
2  C_* \delta \leq \epsilon.
\]
Then, we can choose $\lambda_0 = \lambda_0(\delta, M_1, \omega_0)= \lambda_0(\epsilon, \Lambda, M_0, \omega_0, n) \geq 1$ sufficiently large, such that
\[
 \frac{(M_1 K_\delta +\delta)}{\lambda_0} \leq \delta.
\]
From our choices of $\delta, \lambda_0$, it follows from \eqref{u-v-interior-final} that
\[
\fint_{Q_{3r/2}} |\nabla  u - \nabla v|^2 dxdt \leq 2C_* \delta \leq \epsilon, \quad \forall \ \lambda \geq \lambda_0.
\]
This proves the first assertion of Proposition \ref{interio-approx-lemma}. This assertion, and the fact that $\epsilon$ is sufficiently small imply that
\[
\fint_{Q_{3r/2}} |\nabla v|^2 dx dt \leq \fint_{Q_{3r/4}} |\nabla u|^2 dx dt  + 1 \leq C(n).
\]
From this, and Lemma \ref{L-infty-reg-inter}, it follows that
\[
\norm{\nabla v}_{L^\infty(Q_{r})} \leq C(\Lambda, n) \left( \fint_{Q_{3r/2}} |\nabla v|^2 dx dt\right)^{1/2} \leq C(\Lambda, n),
\]
and hence the second assertion in the proposition is proved, and the proof is complete.
\end{proof}
\subsection{Weighted interior level set estimates}
We begin with the following lemma.
\begin{lemma} \label{density-est-interior-1} Let $\Lambda, M_0, M_1, M_2$ be given positive numbers and let $\beta = \beta(\Lambda, M_0, n) >1$ be as in  \eqref{beta.def}.  Let  $\omega_0: \K \times [0,\infty) \rightarrow (0, \infty)$ be continuous bounded function  satisfying \eqref{omega-limit} and $\norm{\omega_0}_{\infty} \leq M_1$. 
There exists $N = N(\Lambda,  n) >1$ such that for every $\omega\in A_q$, $[\omega]_{A_q} \leq M_2$, with some $1 < q < \infty$,  for every sufficiently small $\epsilon >0$, there exist sufficiently small $\delta = \delta(\epsilon, \Lambda, M_0, M_1, M_2,  n)$  and sufficiently large $\lambda_0 = \lambda_0(\epsilon, \Lambda, M_0, \omega_0,  M_2, n)$ such that the following holds: Suppose that $\A \in \U_{Q_{2R}, \K}( \Lambda, M_0, M_1, \omega_0)$, with its  asymptotical matrix $\tA$ satisfying \eqref{ellip-interior}, for some $R>0$, and for some open interval $\K \subset \mathbb{R}$. Assume that 
\[ [[\A]]_{\textup{BMO}^\#(Q_R, \beta, r_0)} \leq \delta, \quad \text{for some} \quad r_0 \in (0,R),
\] 
then for every $\lambda \geq \lambda_0$, and for weak solution $u \in L^2(\Gamma_{2R}, W^{1,2}(B_{2R}))$  of \eqref{eqn-in}, every $\hat{z} =(\hat{x}, \hat{t}) \in Q_R$, and every $0 < r \leq r_0/5$, if  
\begin{equation} \label{non-empty-iterior-intersection}
Q_r(\hat{z}) \cap  \Big\{ Q_R: \M_{Q_{2R}}(|\nabla u|^2 ) \leq 1\Big \} \cap \Big \{Q_{R}: \M_{Q_{2R}} (|\F|^2) \leq \delta^2 \Big \} \not= \emptyset,
\end{equation}
then 
\[
\omega\Big(\{z \in Q_{R}: \M_{Q_{2R}}(|\nabla u|^2  >N\} \cap Q_r(\hat{z}) \Big) < \epsilon \omega(Q_r(\hat{z})).
\]
\end{lemma}
\begin{proof} The proof of this lemma is standard. However, we present it here for completeness. For a given sufficiently small $\epsilon >0$, choose $\gamma \in (0,1)$ sufficiently small, to be determined, and depending only on $\epsilon$,  $M_2$ and $n$. Let $\delta = \delta(\gamma, \Lambda, M_0, M_1, n) \in (0,1/8)$ and $\lambda_0 =  \lambda_0(\gamma, \Lambda, M_0, \omega_0, n)$ be defined as in Proposition \ref{interio-approx-lemma}. By \eqref{non-empty-iterior-intersection}, there is $z_0 \in Q_r(\hat{z})$ such that
\begin{equation} \label{x-zero-interior}
\M_{Q_{2R}} (|\nabla u|^2)(z_0)  \leq 1, \quad \text{and} \quad \M_{Q_{2R}}(|\F|^2) (z_0) \leq \delta^2.
\end{equation}
Observe that $Q_{3r}(\hat{z}) \subset Q_{4r}(z_0) \subset Q_{2R}$, and
\[
\begin{split}
&  \fint_{Q_{3r}(\hat{z})} |\nabla u|^2  dx dt   \leq \frac{|Q_{4r}(z_0))|}{|Q_{3r}(\hat{z})|}  \fint_{Q_{4r}(z_0)} |\nabla u|^2  dx dt  \leq \left(\frac{4}{3}\right)^{n+2},\\
& \fint_{Q_{3r}(\hat{z})} |\F|^2 dx dt \leq \frac{|Q_{4r}(z_0))|}{|Q_{3r}(\hat{z})|}  \fint_{Q_{4r}(z_0)} |\F|^2 dx dt\leq \left(\frac{4}{3}\right)^{n+2}\delta^2.
\end{split}
\]
These two estimates together with the assumption that $[[\A]]_{\textup{BMO}(Q_1, \beta, r_0)} \leq \delta$  allows us to apply the Proposition \ref{interio-approx-lemma} for $u'(x,t) = u(x,t)/ (4/3)^{n+2}$, and  $\F'(x,t)= \F(x,t)/(4/3)^{n+2}$ to infer that there exists $v \in L^2(Q_{9r/4}(\hat{z}))$ such that
\begin{equation} \label{v-interior-compare-density}
\begin{split}
& \fint_{Q_{9r/4}(\hat{z})}|\nabla u - \nabla v|^2 dxdt  \leq \gamma (4/3)^{n+2}, \quad \text{and} \quad \norm{\nabla w}_{L^{\infty}(Q_{3r/2}(\hat{z}))}  \leq C_*(\Lambda, n).
\end{split}
\end{equation}
Now, let
\[
N = \max\{4C_*(\Lambda, n)^2, 5^{n+2} \}.
\]
We claim that
\begin{equation} \label{set-est-interior}
\begin{split}
& \Big\{Q_r(\hat{z}): \M_{Q_{9r/4}(\hat{z})}(|\nabla (u-v) |^2)  \leq C_*^2\Big\}  \subset \Big\{ Q_r(\hat{z}): \M_{Q_{2R}}(|\nabla u|^2)  \leq N\Big\}.
\end{split}
\end{equation}
To this end, let $z$ to be any point in the set on the left hand side of \eqref{set-est-interior}. We only need to show that
\begin{equation} \label{intereior-N-est}
\M_{Q_{2R}}(|\nabla u|^2)(z)  \leq N.
\end{equation}
We consider the cylinder $Q_\rho(z)$. If $\rho \leq r/2$, we see that $Q_{\rho}(z) \subset Q_{3r/2}(\hat{z}) \subset Q_{2R}$. From this,  it follows
\[
\begin{split}
  \fint_{Q_\rho(z)} \nabla u|^2  dxdt & \leq 2 \left[ \fint_{Q_\rho(z)} |\nabla u - \nabla v|^2  dxdt   +  \fint_{Q_\rho(z)}  |\nabla v|^2   dx dt \right] \\
& \leq 2\Big[\M_{Q_{9r/4}(\hat{z})}(|\nabla u - \nabla v|^2(z)    + \norm{\nabla v}_{L^\infty(Q_{3r/2}(\hat{z}))}^2  \Big]\\
& \leq 2[C_*^2 + C_*^2] = 4C_*^2 \leq N.
\end{split}
\]
On the other hand, if $ \rho > r/2$, we see that $Q_\rho(z) \subset Q_{5\rho}(z_0)$, we can use \eqref{x-zero-interior} to infer that 
\[
 \fint_{Q_\rho(z) \cap Q_{2R}} |\nabla u|^2 dx dt \leq \frac{|Q_{5\rho}(z_0)|}{|Q_\rho(z)|} \fint_{Q_{5\rho}(z_0) \cap Q_{2}} |\nabla u|^2  dxdt \leq 5^{n+2} \leq N.
\]
Hence \eqref{intereior-N-est} follows, which in turns proves \eqref{set-est-interior}. Observe that \eqref{set-est-interior} is equivalent to
\begin{equation} \label{compare-set-interior}
\begin{split}
& \Big\{Q_r(\hat{z}): \M_{Q_{2R}}(|\nabla u|^2) > N\Big\}  \subset E: = \Big\{Q_r(\hat{z}): \M_{Q_{9r/4}(\hat{z})}(|\nabla u-\nabla v|^2) > C_*^2\Big\}.
\end{split}
\end{equation}
Observe that by the weak type (1,1) estimate in Lemma \ref{Hardy-Max-p-p}, and \eqref{v-interior-compare-density}, we see that
\[
|E| \leq C(n) |Q_{9r/4}(\hat{z})| \fint_{Q_{9r/4}(\hat{z})}  |\nabla u - \nabla v|^2  dx dt \leq C_0(n)\gamma |Q_r(\hat{z})|.
\]
Observe that rom Lemma \ref{doubling}, there is 
 $\beta_0 = \beta_0(M_2, n) >0$ such that
 \[
 \omega(E) \leq C(M_2, n) \Big( \frac{|E|}{|Q_r(\hat{z})|} \Big)^{\beta_0} \omega(Q_r(\hat{z})) \leq C^*  \gamma^{\beta_0} \omega(Q_r(\hat{z})),
 \]
where $C^*>0$ is a constant depending only on $M_2$ and $n$.
By choosing $\gamma = \Big(\frac{\e}{C^*}\Big)^{1/\beta_0}$, we obtain the desired result.
\end{proof}
\noindent
We now can estimate a level set of $\M_{Q_{2R}}(|\nabla u|^2)$. This is the main result in this subsection.
\begin{lemma} \label{good-lambda-interior} Let $\Lambda, M_0, M_1, M_2$ be given positive numbers and let $\beta = \beta(\Lambda, M_0, n) >1$ be as in  \eqref{beta.def}.  Let  $\omega_0: \K \times [0,\infty) \rightarrow (0, \infty)$ be a continuous bounded function  satisfying \eqref{omega-limit} and $\norm{\omega_0}_{\infty} \leq M_1$. For every $\omega\in A_q$ with $[\omega]_{A_q} \leq M_2$ for some $1 < q < \infty$,  and for every sufficiently small $\epsilon >0$, let $\delta = \delta(\epsilon, \Lambda, M_0, M_1, M_2,  n)$, $\lambda_0 = \lambda_0(\epsilon, \Lambda, M_0, \omega_0, M_2, n)$, and $N = N(\Lambda,  n) >1$ be defined as in Lemma \ref{density-est-interior-1}. Let $\A \in \U_{Q_{2R}, \K}( \Lambda, M_0, M_1, \omega_0)$ with its  asymptotical matrix $\tA$ satisfying \eqref{ellip-interior} and $[[\A]]_{\textup{BMO}^\#(Q_R, \beta, r_0)} \leq \delta$, and with some open interval $\K \subset \mathbb{R}$ and some $r_0 \in (0,R)$. Then, for every $\lambda \geq \lambda_0$, and every weak solution $u \in L^2(\Gamma_{2R}, W^{1,2}(B_{2R}))$ of \eqref{eqn-in}, if 
\[
\omega \Big (\Big\{Q_R: \M_{Q_{2R}}(|\nabla u|^2) > N \Big\} \Big) \leq \epsilon \omega(Q_R(z)), \quad \forall z \in Q_R, 
\]
it holds that
\[
\begin{split}
 \omega \Big (\Big\{Q_R: \M_{Q_{2R}}(|\nabla u|^2 ) > N \Big\} \Big) &  \leq \epsilon_1  \omega \Big (\Big\{Q_{R}: \M_{Q_{2R}}(|\nabla u|^2 ) > 1 \Big\} \Big)+  \omega \Big (\Big\{Q_R: \M_{Q_{2R}}(|\F|^2) > \delta^2 \Big\} \Big),
\end{split}
\]
where $\e_1 = (20)^{(n+2)q}M_2^2 \e$ defined as in Lemma \ref{Vitali}.
\end{lemma}
\begin{proof} Let us denote
\[
C = \Big\{Q_R: \M_{Q_{2R}}(|\nabla u|^2) > N \Big\},
\]
and
\[
D= \Big\{ Q_R: \M_{Q_{2R}}(|\nabla u|^2 ) >1 \Big\} \cup \Big\{ Q_R : \M_{Q_{2R}}(|\F|^2) > \delta^2 \Big\}.
\]
Clearly, $C\subset D \subset Q_R$. Then in view of Lemma \ref{density-est-interior-1}, we can apply Lemma \ref{Vitali} to obtain the desired estimate.
\end{proof}
\subsection{Proof of Theorem \ref{interior-rg-thm}} By iterating lemma \ref{good-lambda-interior}, we obtain the following result.
\begin{lemma} \label{iteration-interior-lemma} Let $\Lambda, M_0, M_1, M_2$ be given positive numbers and let $\beta = \beta(\Lambda, M_0, n) >1$ be as in  \eqref{beta.def}.  Let  $\omega_0: \K \times [0,\infty) \rightarrow (0, \infty)$ be a continuous bounded function  satisfying \eqref{omega-limit} and $\norm{\omega_0}_{\infty} \leq M_1$. Let $\omega\in A_q$ with $[\omega]_{A_q} \leq M_2$ for some $1 < q < \infty$. For every sufficiently small $\epsilon >0$, let $\delta = \delta(\epsilon, \Lambda, M_0, M_1, M_2,  n)$, $\lambda_0 = \lambda_0(\epsilon, \Lambda, M_0, \omega_0,  M_2, n)$, and $N = N(\Lambda,  n) >1$ be defined as in Lemma \ref{density-est-interior-1}. Suppose that $\A \in \U_{Q_{2R}, \K}( \Lambda, M_0, M_1, \omega_0)$ with its  asymptotical matrix $\tA$ satisfying \eqref{ellip-interior} and $[[\A]]_{\textup{BMO}^\#(Q_R, \beta, r_0)} \leq \delta$, and with some open interval $\K \subset \mathbb{R}$ and some $r_0 \in (0,R)$. Then, if $u \in L^2(\Gamma_{2R}, W^{1,2}(B_{2R}))$ is a weak solution of \eqref{eqn-in} such that 
\[
\omega \Big( \Big\{Q_R: \M_{Q_{2R}} (|\nabla u|^2) > N \Big\} \Big) \leq \epsilon \omega(Q_{R}(z)), \quad \forall \ z \in Q_R,
\]
then for $\epsilon_1 = (20)^{(n+2)q}M_2^2\epsilon$, it holds that
\[
\begin{split}
\omega\Big (\Big\{Q_R: \M_{Q_{2R}} (|\nabla u|^2) > N^k \Big\} \Big) & \leq \epsilon_1^k  \omega\Big(\Big\{Q_R: \M_{Q_{2R}}(|\nabla u|^2) > 1 \Big\} \Big)  +  \sum_{i=1}^k \epsilon_1^i \omega \Big(\Big\{Q_R: \M_{Q_{2R}}(|\F|^2) > \delta^2 N^{k-i} \Big\} \Big) .
\end{split}
\]
\end{lemma}
\begin{proof} We skip the proof because it is the same as that of Lemma \ref{global-iterating-lemma} below.
\end{proof} \noindent
After the accomplishment of Lemma \ref{iteration-interior-lemma}, the rest of the proof of Theorem \ref{interior-rg-thm} is similar to that of Theorem \ref{main-theorem} for boundary regularity estimates in the next section.  We therefore skip it.
\section{Weighted $W^{1,p}$-regularity estimates on flat domains} \label{global-regularity-section}
For each $r >0$ and for $x_0~=~(x^{0}_1, x^{0}_2, \cdots, x^{0}_n)~\in~\mathbb{R}^n$, let us denote
\[
\begin{split}
B_r^+(x_0)  = \{ y = (y_1, y_2, \cdots, y_n) \in B_r(x_0):\ y_n >x^0_{n}\}, \quad B_r^+ = B_r^+(0), \\
T_r (x_0) = \{ x =  (x_1, x_2, \cdots, x_n) \in \partial B_{r}^+(0): \ x_n =x_{n}^{0}\}, \quad T_r = T_r(0).
\end{split}
\]
For given $z_0 = (x_0, t_0) \in \R^{n+1}$ and for given $r >0$, we denote
\[
Q_r^+(z_0) = B_r(x_0)^+ \times \Gamma_r(t_0), \quad \text{where} \quad \Gamma_r(t_0) = (t_0 -r^2, t_0].
\]
Similarly as before, when $z_0 =0$, we write
\[
Q_r^+ = Q_r^+(0), \quad \Gamma_r = \Gamma_r(0).
\]
In this section, for $R >0, \lambda >0$, $\A \in \U_{Q_{2R}^+, \K}(\Lambda, M_0, M_1, \omega_0)$, with some given $\omega_0$, we study weak solution $u$ of the problem
\begin{equation} \label{boundary-eqn-r}
\left\{
\begin{array}{cccl}
u_t - \text{div}[\A_\lambda(x, t, u, \nabla u)]  & = & \text{div}[{\bf F}] & \quad \text{in} \quad \ Q_{2R}^+, \\
u & =& 0 & \quad \text{on} \quad T_{2R} \times \Gamma_{2R},
\end{array}
\right.
\end{equation}
\noindent
where $\A_\lambda$ is defined in \eqref{A-lambda}. Observe that the equation in \eqref{boundary-eqn-r} is reduced to \eqref{main-eqn} when taking $\lambda =1$. As we already discussed, enlarging the class of equations in \eqref{main-eqn} to the class of equations in \eqref{boundary-eqn-r} with the parameter $\lambda>0$ is essential in our approach. This is mainly because of the homogeneity in the Calder\'on-Zygmund type regularity estimates.

Let $\tA$ be the asymptotical matrix of $\A$, we recall that $\tA = A +D$, where $A$ is an $n\times n $ symmetric matrix, and $D$ is a $n\times n$ skew-symmetric matrix. Moreover, $A \in L^\infty(Q_{2R}^+)$,  $D \in L^\infty(\Gamma_{2R}, \textup{BMO})$, and
\begin{equation} \label{ellip-bd}
\Lambda^{-1} |\eta|^2  \leq \wei{A(x,t) \eta, \eta}, \quad \norm{A}_{L^\infty(Q_{2R}^+)} \leq \Lambda, \quad  
\norm{D}_{L^\infty(\Gamma_r, \textup{BMO})} \leq M_0, \quad \text{for a.e.} \ (x,t) \in Q_{2R}^+, \quad \forall \ \eta \in \R^n.
\end{equation}
In this section, let  $\gamma = \gamma(\Lambda, M_0, n)>0$ be defined as in  Lemma \ref{Gehring-2}, and let 
\begin{equation} \label{beta-2.def}
\beta = \frac{2(2+\gamma)}{\gamma} > 1.
\end{equation}
For some $r_0 \in (0,R)$,  we denote
\[
[[\tA]]_{\textup{BMO}^\#(Q_R^+, \beta, r_0)} = \left( \sup_{z_0 = (x_0, t_0) \in Q_R^+} \sup_{0 < \rho < r_0} \fint_{Q_{\rho}(z_0) \cap Q_{2R}^+} |\tA(x,t) -\wei{\tA}_{B_\rho(x_0) \cap B_{2R}^+}(t)|^\beta dxdt \right)^{1/\beta},
\]
where $\wei{\tA}_{B_\rho(x_0) \cap B_{2R}^+}(t)$ denotes the mean of $\tA$ on $B_\rho(x_0)\cap B_{2R}^+$. Our main result of the section is 
the following theorem, which in turns also gives Theorem \ref{W-1p-reg-b} when taking $\lambda =1$.
\begin{theorem} \label{main-theorem}   Let $\Lambda >0 , M_0 >0 , M_1 >0, M_2 \geq 1$ and $p > 2$. Then there are $\delta = \delta (\Lambda, M_0, M_1, M_2, p, n)$ sufficiently small and $\beta = \beta(\Lambda, M_0, n) >1$ such that  the following holds: Suppose that $R>0$, $\K \subset \R$ is an open interval and $\omega_0:\K\times [0, \infty)\rightarrow [0~,~\infty)$ is continuous satisfying \eqref{omega-limit} and $\norm{\omega_0}_{\infty} \leq M_1$.  Suppose also that $\A \in \mathbb{U}_{Q_{2R}^+, \mathbb{K}} (\Lambda, M_0, M_1, \omega_0)$, with its asymptotical matrix $\tA$ satisfying \eqref{ellip-bd} and 
\begin{equation}  \label{smallness-BMO}
[[\tA]]_{\textup{BMO}\#(Q_{2R}^+, \beta, R)} \leq \delta.
\end{equation}
If $u \in L^2(\Gamma_{2R}, W^{1,2}(B_{2R}^+))$ is a weak solution of  \eqref{boundary-eqn-r} with some $\lambda >0$, then the estimate
\begin{equation} \label{W-1p-flat}
\begin{split}
\left( \frac{1}{\omega(Q_R^+)}\int_{Q_{R}^+}   |\nabla u|^p \omega(x,t) dx dt \right)^{1/p}& \leq C \left[ \left( \frac{1}{\omega(Q_{2R}^+)} \int_{Q_{2R}^+} | \F|^p \omega(x, t) dx dt \right)^{1/p}\right. \\
& \quad \quad + \left.  \max \left\{ \left( \fint_{Q_{2R}^+} |\nabla u|^2 dx dt \right)^{1/2}, \lambda^{-1} \right\} \right]
\end{split}
\end{equation}
holds if $|\F| \in L^2(Q_{2R}^+) \cap L^p(Q_{2R}^+, \omega)$, and for some uniform constant $C$ depending only on $\Lambda, M_0, M_1, M_2, \omega_0, n, p$, and for some $\omega \in A_{p/2}$ with $[\omega]_{A_{p/2}} \leq M_2$.
\end{theorem}

\subsection{Boundary approximation estimates} \label{boundary-global-section} The following Proposition is the main result of this subsection.
\begin{proposition} \label{boundary-approximation-lemma} 
Let $\Lambda >0, M_0 > 0, M_1 >0$ be given. Let $\beta$ be as in \eqref{beta-2.def}, and let  $\omega_0$ be  a continuous function satisfying \eqref{omega-limit}, $\norm{\omega_0}_{\infty} \leq M_1$. Then, for every small number $\epsilon >0$, there exist $\delta = \delta(\epsilon, \Lambda, M_0, M_1, n) >0$ sufficiently small and sufficiently large $\lambda_0 =  \lambda_0(\epsilon, \Lambda, M_0, \omega_0, n) \geq 1$  such that the following holds:  Assume that $\A \in \U_{Q_{2R}^+, \K}( \Lambda, M_0, M_1, \omega_0)$, with its  asymptotical matrix $\tA$ satisfying \eqref{ellip-bd} for some $R>0$ and some open interval $\K \subset \mathbb{R}$. Assume also that
\[
[[\tA]]_{\textup{BMO}^\#(Q_R^+, \beta, r_0)} \leq \delta, \quad \text{for some} \quad r_0 \in (0,R).
\]
Then, for some $r \in (0,r_0/2)$ some $z_0 = (x_0, t_0) \in T_R \times \Gamma_R$ and for $\lambda \geq \lambda_0$, if  
\[
\fint_{Q_{2r}^+(z_0)} | \F(x,t)|^2 dx dt  \leq \delta,
\]
and  if  $u$ is a weak solution of \eqref{boundary-eqn-r} satisfying
\[
\fint_{Q_{2r}^+(z_0)} |\nabla u|^2 dx dt \leq 1,
\]
then there is $v \in L^2(Q_{3r/2}^+(z_0))$, and constant $C = C(\Lambda, n)$ such that
\[
\fint_{Q_{3r/2}^+(z_0)} |\nabla u - \nabla v|^2 dx dt \leq \epsilon,
\quad \text{and} \quad \norm{\nabla v}_{L^\infty(Q_{r}^+(z_0))} \leq C.
\]
\end{proposition}
\begin{proof} The proof is similar to that of Proposition \ref{interio-approx-lemma} using two steps of approximation. Essential ingredients are Theorem \ref{existence-theorem}, Lemma \ref{Gehring-2}, and Lemma \ref{L-infty-reg-bdr}. We skip the proof.
\end{proof}
\noindent
\subsection{Weighted boundary level set estimates}

We need the following lemma which is a restated version of Lemma \ref{density-est-interior-1} for interior cylinders in $Q_{2R}^+$.

\begin{lemma} \label{density-est-bdry-1} Let $\Lambda, M_0, M_1, M_2$ be given positive numbers and let $\beta = \beta(\Lambda, M_0, n) >1$ be as in  \eqref{beta.def}.  Let  $\omega_0: \K \times [0,\infty) \rightarrow (0, \infty)$ continuous bounded function  satisfying \eqref{omega-limit} and $\norm{\omega_0}_{\infty} \leq M_1$. 
There exists $N_1 = N_1(\Lambda,  n) >1$ such that for every $\omega\in A_q$ with $[\omega]_{A_q} \leq M_2$ for some $1 < q < \infty$,  for every sufficiently small $\epsilon >0$, there exist sufficiently small $\delta_1 = \delta_1(\epsilon, \Lambda, M_0, M_1, M_2,  n)$  and sufficiently large $\lambda_1 = \lambda_1(\epsilon, \Lambda, M_0, \omega_0,  M_2, n)$ such that the following holds: Suppose that $\A \in \U_{Q_{2R}^+, \K}( \Lambda, M_0, M_1, \omega_0)$, with its  asymptotical matrix $\tA$ satisfying \eqref{ellip-bd}, for some $R>0$, and for some open interval $\K \subset \mathbb{R}$. Assume that 
\[ [[\A]]_{\textup{BMO}^\#(Q_R^+, \beta, r_0)} \leq \delta, \quad \text{for some} \quad r_0 \in (0,R), \quad\text{and some} \quad \delta \leq \delta_1,
\] 
Assume also that  $u \in L^2(\Gamma_{2R}, W^{1,2}(Q_{2R}^+)) $ is a weak solution of \eqref{boundary-eqn-r} with $\lambda \geq \lambda_1$, and assume that 
\begin{equation} \label{non-empty-bdry-intersection-1}
Q_r(\hat{z}) \cap  \Big\{ Q_R^+: \M_{Q_{2R}^+}(|\nabla u|^2 ) \leq 1\Big \} \cap \Big \{Q_{R}^+: \M_{Q_{2R}^+} (|\F|^2) \leq \delta^2 \Big \} \not= \emptyset,
\end{equation}
for some $0 < r \leq r_0/5$ and for some $\hat{z} \in Q_R^+$ with  $Q_{3r}(\hat{z}) \subset Q_{2R}^+$. Then 
\[
\omega\Big(\{x \in Q_{R}^+: \M_{Q_{2R}^+}(|\nabla u|^2)  >N_1\} \cap Q_r(\hat{z}) \Big )< \epsilon \omega (Q_r(\hat{z})).
\]
\end{lemma}
\noindent
Our next lemma is similar to Lemma \ref{density-est-bdry-1}, but  for cylinders centered on the flat boundary $T_1 \times \Gamma_1$ of $Q_2^+$.
\begin{lemma} \label{density-est-bdry-2} Let $\Lambda, M_0, M_1, M_2$ be positive number and let $\beta, \omega_0$ be as in Proposition \ref{boundary-approximation-lemma}. Also, let $\epsilon >0$ sufficiently small, and $\omega \in A_q$ with $[\omega]_{A_q} \leq M_2$ for some $1 < q < \infty$. There exist $N_2 = N_2(\Lambda, n) >1$, a sufficiently small number $\delta_2 = \delta_2(\epsilon, \Lambda, M_0, M_1, M_2, n)>0$, and a sufficiently large number $\lambda_2 = \lambda_2(\epsilon, \Lambda, M_0, \omega_0,  M_2, n)$ such that  the following statement holds true. 
Suppose that $\A \in \U_{Q_{2R}^+, \K}( \Lambda, M_0, M_1, \omega_0)$, with its  asymptotical matrix $\tA$ satisfying \eqref{ellip-bd}, for some $R>0$, and for some open interval $\K \subset \mathbb{R}$. Assume that 
\[ [[\A]]_{\textup{BMO}^\#(Q_R^+, \beta, r_0)} \leq \delta, \quad \text{for some} \quad r_0 \in (0,R), \quad\text{and some} \quad \delta \leq \delta_2,
\] 
Assume also that  $u \in L^2(\Gamma_{2R}, W^{1,2}(Q_{2R}^+)) $ is a weak solution of \eqref{boundary-eqn-r} with $\lambda \geq \lambda_2$
and
\begin{equation} \label{non-empty-bdry-intersection-2}
Q_r(z_0) \cap  \Big\{ Q_R^+: \M_{Q_{2R}^+}(|\nabla u|^2  ) \leq 1\Big \} \cap \Big \{Q_{R}^+: \M_{Q_{2R}^+} (|\F|^2) \leq \delta^2 \Big \} \not= \emptyset,
\end{equation}
for some $0 < r \leq r_0/5$ and for some $z_0 \in T_1\times \Gamma_1$. Then 
\[
\omega \Big(\{x \in Q_{R}^+: \M_{Q_{2R}^+}(|\nabla u|^2)  >N_2\} \cap Q_r(z_0) \Big) < \epsilon \omega(Q_r(z_0)).
\]
\end{lemma}
\begin{proof} The proof is similar to that of Lemma \ref{density-est-interior-1} using Proposition \ref{boundary-approximation-lemma} instead of Proposition \ref{interio-approx-lemma}. We therefore skip the proof.
\end{proof}
Now, combining the previous Lemma, Lemma \ref{density-est-bdry-1} and Lemma \ref{density-est-bdry-2} we can prove the following result, which is also the main ingredient for the estimates of our level sets.
\begin{lemma} \label{density-est-bdry-3}  Let $\Lambda, M_0, M_1, M_2$ be given positive numbers an let $\omega_0: \K \times [0,\infty) \rightarrow (0, \infty)$ continuous bounded function  satisfying \eqref{omega-limit} and $\norm{\omega_0}_{\infty} \leq M_1$.  Then, for $\epsilon >0$ sufficiently small, for every $\omega \in A_q$ with $[\omega]_{A_q} \leq M_2$ for some $1 < q < \infty$, there are  $\beta = \beta(\Lambda, M_0,n) >1$, $N = N(\Lambda, n) >1$,  sufficiently small positive number $\delta = \delta(\epsilon, \Lambda, M_1, M_0, M_2, q, n)$, and sufficiently large number $\lambda_0 =\lambda_0(\epsilon, \Lambda, M_0, \omega_0, M_2, n)$ such that  the following statement holds true. Suppose that $\A \in \U_{Q_{2R}^+, \K}( \Lambda, M_0, M_1, \omega_0)$, with its  asymptotical matrix $\tA$ satisfying \eqref{ellip-bd}, for some $R>0$, and for some open interval $\K \subset \mathbb{R}$. Assume that 
\[ [[\A]]_{\textup{BMO}^\#(Q_R^+, \beta, r_0)} \leq \delta, \quad \text{for some} \quad r_0 \in (0,R).
\] 
Assume also that  $u \in L^2(\Gamma_2, W^{1,2}(Q_2^+)$ is a weak solution of \eqref{boundary-eqn-r} for $\lambda \geq \lambda_0$. If  $\hat{z} \in \overline{Q}_1^+$ and $r \in (0, r_0/20)$ such that
\begin{equation} \label{bdry-intersection-3}
\omega \Big(  \Big\{ Q_R^+: \M_{Q_{2R}^+}(|\nabla u|^2   ) > N \Big \} \cap 
Q_r(\hat{z}) \Big ) \geq \epsilon \omega(Q_r(\hat{z})),
\end{equation}
then 
\[
Q_r(\hat{z}) \cap Q_R^+ \subset \Big\{ Q_1^+: \M_{Q_{2R}^+}(|\nabla u|^2 ) > 1 \Big \} \cup \Big\{Q_1^+: \M_{Q_{2R}^+}(|\F|^2) > \delta^2 \Big\}.
\]
\end{lemma}
\begin{proof} Let $\e' =\frac{\epsilon}{M_2 4^{(n+2) q}}$. Let $N = \max\{N_1, N_2\}$, where $N_1, N_2$ are respectively defined in Lemma \ref{density-est-bdry-1} and Lemma \ref{density-est-bdry-2}. Also, let
\[
\begin{split}
\delta & = \min\{\delta_1(\epsilon', \Lambda, M_0, M_1, M_2, n), \delta_2(\e, \Lambda, M_0,M_1, M_2, n)\},\\
\lambda_0 & = \max\{\lambda_1(\epsilon', \Lambda, M_0, \omega_0, M_2, n), \lambda_2(\epsilon, \Lambda, M_0, \omega_0,  M_2, n) \},
\end{split}
\]
where $\delta_1, \delta_2, \lambda_1, \lambda_2$ are again defined in Lemma \ref{density-est-bdry-1} and Lemma \ref{density-est-bdry-2} respectively. We prove Proposition \ref{density-est-bdry-3} holds with our choice of $N, \delta, \lambda_0$. Let us denote $\hat{z} = (\hat{x}, \hat{t})$. Observe that if $B_{3r}(\hat{x}) \subset B_{2R}^+$, then the conclusion of our proposition follows directly from Lemma \ref{density-est-bdry-1}. Therefore, it remains to consider the case that $B_{3r}(\hat{x}) \cap T_R \not= \emptyset$. In this case, we write $\hat{x} = (\hat{x}', \hat{x}_n)$ and then let $\hat{x}_0 = (\hat{x}', 0)\in B_{3r}(\hat{x}) \cap T_R$. We assume by contradiction that there is $z_0 = (x_0, t_0) \in Q_r(\hat{z})\cap Q_R^+$ such that
\begin{equation} \label{x-zero-contra}
\M_{Q_{2R}^+}(|\nabla u|^2  ) (z_0) \leq 1,\quad \text{and} \quad \M_{Q_{2R}^+}(|\F|^2)(z_0) \leq \delta^2.
\end{equation}
Observe that for $\rho = 4r$,
\[
z_0 \in Q_r(\hat{z}) \cap Q_R^+ \subset Q_{\rho}^+(\hat{z}_0) \cap Q_R^+, \quad \text{where} \quad \hat{z}_0 = (\hat{x}_0, \hat{t}) \in T_R \times \Gamma_R.
\]
This and \eqref{x-zero-contra} particular imply that
\[
Q_\rho^+(\hat{z}_0) \cap  \Big\{Q_R^+ : \M_{Q_{2R}^+}(|\nabla u|^2 ) \leq 1 \Big\} \cap \Big\{Q_R^+: \M_{Q_{2R}^+}(|\F|^2) \leq \delta^2\Big\} \not=\emptyset.
\]
Moreover, since $r < r_0/20$, we see that $\rho< r_0/5$. Hence, from our choice of $N, \delta$, we can apply Lemma \ref{density-est-bdry-2} to conclude that
\[
\begin{split}
\omega \Big(\{ Q_{R}^+: \M_{Q_{2R}}(|\nabla u|^2)  >N \} \cap Q_r(\hat{z}) \Big)
& \leq \omega \Big(\{ Q_{R}^+: \M_{Q_{2R}}(|\nabla u|^2+ |\pi|^2 )  >N_2\} \cap Q_\rho(\hat{z}_0) \Big) \\
& <  \epsilon' \omega (Q_\rho(\hat{z}_0)) \leq \epsilon \omega(Q_r(\hat{z})),
\end{split}
\]
where we have used Lemma \ref{doubling} in the last estimate. Note that the last estimate contradicts to \eqref{bdry-intersection-3}. The proof of the proposition is then complete.
\end{proof}
Lemma \ref{density-est-bdry-3} implies the following important result, which is the main result of the subsection.
\begin{lemma} \label{good-lambda-bdr} Let $\Lambda, M_0, M_1, M_2$ be positive numbers and let $\beta, \omega_0$ be as in Lemma \ref{density-est-bdry-3}. Also, let $\omega \in A_q$ with some $1 < q < \infty$ with $[\omega]_{A_q} \leq M_2$. For any $\epsilon >0$, let $\delta, N, \lambda_0$ be defined as in Lemma \ref{density-est-bdry-3}.  Suppose that $\A \in \U_{Q_{2R}^+, \K}( \Lambda, M_0, M_1, \omega_0)$, with its  asymptotical matrix $\tA$ satisfying \eqref{ellip-bd}, for some $R>0$, and for some open interval $\K \subset \mathbb{R}$. Assume that 
\[ [[\A]]_{\textup{BMO}^\#(Q_R^+, \beta, r_0)} \leq \delta, \quad \text{for some} \quad r_0 \in (0,R).
\] 
Assume also that  $u \in L^2(\Gamma_2, W^{1,2}(Q_{2R}^+)$ is a weak solution of \eqref{boundary-eqn-r} for $\lambda \geq \lambda_0$, and 
\[
\omega \Big(\Big\{z \in Q_R^+: \M_{Q_{2R}^+}(|\nabla u|^2) (z) > N \Big\} \Big) \leq \epsilon \omega (Q_R(z)), \quad \forall \ z \in Q_R^+,
\]
then for $\epsilon_1 = (20)^{(n+2)q} M_2^2 \epsilon$,
\[
\begin{split}
& \omega \Big(\Big\{Q_R^+: \M_{Q_{2R}^+}(|\nabla u|^2 > N \Big\} \Big) \\
 & \leq \epsilon_1\omega \Big( \Big\{Q_R^+: \M_{Q_{2}^+}(|\nabla u|^2 ) > 1 \Big\} \Big) + \omega \Big(\Big\{Q_R^+: \M_{Q_{2R}^+}(|\F|^2) > \delta^2 \Big\} \Big) .
\end{split}
\]
\end{lemma}
\begin{proof} The same as that of Lemma \ref{good-lambda-interior}, using 
Lemma \ref{density-est-bdry-3} and the modified Vitali's covering lemma, Lemma \ref{Vitali}.
\end{proof}
\begin{proof}[Proof of Theorem \ref{main-theorem}: Weighted $W^{1,p}$-regularity estimate on flat domain] From Lemma \ref{good-lambda-bdr}, Theorem \ref{main-theorem} is generally expected to follow. However, due to our new ingredient related to the parameter $\lambda$, a careful scaling argument has to be performed, and details therefore needed. By iterating Lemma \ref{good-lambda-bdr}, we obtain the following lemma
\begin{lemma} \label{global-iterating-lemma} Let $A, \Lambda, M_{0}, M_1, M_2$ be positive, $q > 1$, and  let $\epsilon >0$ sufficiently small. Also, let  $\omega_0: \K \times [0, \infty) \rightarrow [0, \infty)$ be continuous satisfing $\norm{\omega_0}_{\infty} \leq M_1$, and let $N = N(\Lambda, n)$,  $\delta = \delta (\epsilon,  \Lambda, M_0, M_1, M_2,  n)$, and $\lambda_0 = \lambda_0(\epsilon, \Lambda, M_0,  \omega_0, M_2, n)$ be as in Lemma \ref{good-lambda-bdr}.  For some $R>0$,  and for some open interval $\K \subset \R$, suppose $\A \in \U_{Q_{2R}^+, \K}(\Lambda, M_0, M_1, \omega_0)$ with its asymptotical matrix $\tA$  satisfying 
\eqref{ellip-bd} and 
\[ [[\A]]_{\textup{BMO}^\#(Q_R^+, \beta, r_0)} \leq \delta, \quad \text{for some} \quad r_0 \in (0,R).
\] 
Assume also that  $u \in L^2(\Gamma_2, W^{1,2}(Q_{2R}^+)$ is a weak solution of \eqref{boundary-eqn-r} for $\lambda \geq \lambda_0$ satisfying
\begin{equation}\label{bdry-cond-last}
\omega(\{Q_R^+: \M_{ Q_{2R}^+}(|\nabla u|^{2})  > N \}) \leq \epsilon \omega(Q_{R}^+(z)), \quad \forall \ z \in \overline{Q}_R^+,
\end{equation}
for some $\omega \in A_q$ with $[\omega]_{A_q} \leq M_2$. Then for $\epsilon_1 = 20^{(n+2)q}  M_2^2\epsilon$, and for any $k \in \mathbb{N}$, 
\begin{equation} \label{iteration-formula}
\begin{split}
 \omega \Big  (\Big\{Q_R^+: \M_{Q_{2R}^+}(|\nabla u|^2 > N^k \Big\} \Big) & \leq \epsilon_1^k \omega\Big (\Big\{Q_R^+: \M_{Q_{2R}^+}(|\nabla u|^2 )> 1 \Big\} \Big) \\
 & \quad  + \sum_{i=1}^k \e_1^i \omega\Big (\Big\{Q_{R}^+: \M_{Q_{2R}^+}(|\F|^2) > \delta N^{k-i}\Big\} \Big) .
\end{split}
\end{equation}
\end{lemma}
\begin{proof} We use induction on $k$. If $k =1$, \eqref{iteration-formula}  holds as a result of Lemma \ref{good-lambda-bdr}. Now, let us assume that Lemma \ref{global-iterating-lemma} holds for some $k \in \{1, 2,\cdots, k_0\}$ with some $k_0 \in \N$. Assume that $u$ is a weak solution of \eqref{boundary-eqn-r} with some $\lambda \geq \lambda_0$ so that \eqref{bdry-cond-last} holds. 
Now, let us define $u' = u/\sqrt{N}$, $\F' = \F/\sqrt{N}$, and $\lambda' = \lambda \sqrt{N} \geq \lambda \geq \lambda_0$. Then, we see that $u'$ is a weak solution of
\[
\left\{
\begin{array}{cccl}
\text{div}[\A_{\lambda'}(x,u', \nabla u')] & = & \text{div}[\F'], & \quad \text{in} \quad Q_{2R}^+,\\
u' & =& 0, & \quad \text{on} \quad T_{2R} \times \Gamma_{2R}.
\end{array} \right.
\]
Moreover, $\forall \ z \in \overline{Q}_{R}^+$
\[
\omega\Big(\Big\{Q_R^+: \M_{Q_{2R}^+}(|\nabla u'|^2 ) > N  \Big\} \Big) = \omega\Big(\Big\{Q_R^+: \M_{Q_{2R}^+}(|\nabla u|^2 ) > N^2  \Big\} \Big) \leq \epsilon \omega(Q_{R}^+(z)).
\]
Therefore, we can apply the induction hypothesis to conclude that
\[
\begin{split}
\omega \Big(\Big\{Q_R^+: \M_{Q_{2R}^+} (|\nabla u'|^2 ) > N^l \Big\} \Big) & \leq \epsilon_1^l \omega \Big(\Big\{Q_R^+: \M_{Q_{2R}^+}(|\nabla u'|^2 ) > 1 \Big\} \Big) \\
 & \quad +  \sum_{i=1}^l \epsilon_1^i \omega \Big(\Big\{Q_R^+: \M_{Q_{2R}^+}(|\F'|^2) > \delta^2 N^{l-i} \Big\} \Big).
\end{split}
\]
From this, and by changing $u'$ back to $u$ and using the case $k =1$, we obtain \eqref{iteration-formula} for $k = l+1$. The proof is therefore complete.
\end{proof}
\noindent
We now complete the proof. Let $N = N(\Lambda, n)$ and $\beta = \beta(\Lambda, M_0, n)$ be defined as in Lemma \ref{global-iterating-lemma}. For $p >2$, we denote $q = p/2 >1$, and choose $\epsilon >0$ and sufficiently small and depending only on $\Lambda, M_2, n, p$ such that
\begin{equation} \label{choice-epsi}
\epsilon_1 N^{q} = 1/2,
\end{equation}
where $\epsilon_1$ is defined in Lemma \ref{global-iterating-lemma}. With this $\epsilon$, we can now choose 
\[ 
\begin{split}
\delta = \delta(\Lambda, M_0, M_1, M_2, p, n), \quad \lambda_0 = \lambda_0(\Lambda,  M_0,  \omega_0, M_2, p, n)
 \end{split}
\]
as determined by Lemma \ref{global-iterating-lemma}.  Now, assume \eqref{smallness-BMO} holds with this choice of $\delta$.  Let $u $ be a weak solution of \eqref{boundary-eqn-r} with $\lambda >0$.
We first prove the estimate in Theorem \ref{main-theorem} with the extra condition that
\begin{equation} \label{extra}
\omega\Big(\Big\{Q_R^+: \M_{Q_{2R}^+} (|\nabla u|^2) > N \Big\} \Big) \leq \epsilon \omega (Q_{R}(z)), \quad \forall \ z \in Q_R^+,
\end{equation}
\noindent
and then show how to remove this condition at the end. To perform the  proof with \eqref{extra}, we consider two  cases depending on whether $\lambda \geq \lambda_0$ or not. \\ \ \\
{\bf Step I:} We assume now that $\lambda \geq \lambda_0$, and assume also that \eqref{extra} holds. Let us consider the sum
\[
S = \sum_{k=1}^\infty N^{q k}  \omega\Big (\Big\{Q_R^+ : \M_{Q_{2R}^+}(|\nabla u|^2) > N^k \Big\} \Big).
\]
From \eqref{extra}, we can apply Lemma \ref{global-iterating-lemma} to obtain
\[
\begin{split}
S & \leq \sum_{k=1}^\infty N^{kq}  \sum_{i=1}^k \epsilon_1^i \omega\Big( \Big\{ Q_R^+: \M_{Q_{2R}^+}(|\F|^2) > \delta N^{k-i} \Big\}\Big)  + \sum_{k=1}^\infty \Big(N^{q} \epsilon_1\Big)^k \omega \Big( \Big\{Q_R^+: \M_{Q_{2R}^+}(|\nabla u|^2 ) >1 \Big\}\Big).
\end{split}
\]
By Fubini's theorem, the above estimate can be rewritten as
\begin{equation} \label{Fubini-express}
\begin{split}
S &\leq \sum_{j=1}^\infty (N^q \epsilon_1)^j  \sum_{k=j}^\infty N^{q(k-j)} \omega \Big(\Big\{ Q_R^+: \M_{Q_{2R}^+}(|\F|^2) > \delta N^{k-j} \Big\}\Big) \\
& \quad + \sum_{k=1}^\infty \Big(N^{q} \epsilon_1\Big)^k \omega\Big( \Big\{Q_R^+: \M_{Q_{2R}^+}(|\nabla u|^2) >1 \Big\}\Big).
\end{split}
\end{equation}
Observe that 
\[
\omega \Big (\Big\{Q_R^+: \M_{Q_{2R}^+}(|\nabla u|^2 ) >1 \Big\}\Big) \leq\omega (Q_R^+).
\]
From this, the choice of $\epsilon$ as in \eqref{choice-epsi}, Lemma \ref{measuretheory-lp}, and \eqref{Fubini-express} it follows that
\[
\begin{split}
S \leq C \left [ \norm{\M_{Q_{2R}^+} ( |\F|^2)}_{L^q(Q_R^+, \omega)}^q + \omega(Q_R^+) \right].
\end{split}
\]
Applying Lemma \ref{Hardy-Max-p-p}, and Lemma \ref{measuretheory-lp} again, we see that
\[
\norm{\M_{Q_{2R}^+}(|\nabla u|^2)}_{L^q(Q_R^+, \omega)}^q \leq C\left [  \norm{\F}_{L^p(Q_{2R}^+, \omega)}^p + \omega(Q_R^+) \  \right].
\]
By the Lesbegue's differentiation theorem, we observe that
\[
|\nabla u(x,t)|^2  \leq \M_{Q_{2R}^+}(|\nabla u|^2)(x,t), \quad \text{a.e} \ (x,t)\ \in Q_R^+.
\]
Hence,
\begin{equation} \label{L-p-1}
\norm{\nabla u}_{L^p(Q_R^+, \omega)}^p \leq C\left[ \norm{\F}_{L^p(Q_{2R}^+, \omega)}^p+ \omega(Q_R^+)\right].
\end{equation}
Thus, we have proved Theorem \ref{main-theorem} as long as $u$ is a weak solution of \eqref{boundary-eqn-r} with $\lambda \geq \lambda_0$ and \eqref{extra} holds. \\ \ \\
{\bf Step II:} We study the case $0 < \lambda < \lambda_0$. Assume also that  $u$ is a weak solution of \eqref{boundary-eqn-r} and \eqref{extra} holds.  Let us denote $u' = u/(\lambda_0/\lambda), \F' = \F/(\lambda_0/\lambda)$. Then, $u'$ is a weak solution of
\[
\left\{
\begin{array}{cccl}
\text{div}[\A_{\lambda_0}(x, u', \nabla u')] & = & \text{div}[\F'], & \quad \text{in} \quad Q_{2R}^+, \\
u' & = & 0, & \quad \text{on} \quad T_{2R} \times \Gamma_{2R}.
\end{array}
\right.
\]
Moreover, because of \eqref{extra} and $\lambda_0/\lambda \geq 1$, we also have
\[
\omega\Big(Q_R^+ : \Big\{ \M_{Q_{2R}^+} (|\nabla u'|^2) >N \Big\}\Big) \leq \epsilon \omega(Q_{R}(z)), \quad \forall \ z \in \overline{Q}_R^+.
\]
Therefore, applying the conclusion of \eqref{L-p-1} for $u'$, we also obtain
\[
\norm{\nabla u'}_{L^p(Q_R^+, \omega)} \leq C\left[\norm{\F'}_{L^p(Q_{2R}, \omega)} + \omega(Q_R)^{1/p}  \right].
\]
Thus, 
\begin{equation*} 
\norm{\nabla u}_{L^p(Q_R^+, \omega)} \leq C\left[\norm{\F/\mu}_{L^p(Q_{2R}^+, \omega)}  + \lambda_0\omega(Q_R)^{1/p}/\lambda   \right].
\end{equation*}
\noindent
{\bf Final step:} Up to now from the first two steps, we have proved that if $u$ is a weak solution of \eqref{boundary-eqn-r} with $\lambda>0$ and if \eqref{extra} holds, then 
\begin{equation} \label{ineq-extra-cond-2}
\norm{\nabla u}_{L^p(Q_R^+, \omega)} \leq C\left[\norm{\F}_{L^p(Q_{2R}^+, \omega)}  + \omega(Q_R)^{1/p}\max\{\lambda^{-1}, 1\}   \right].
\end{equation}
It therefore only remains to remove the extra condition \eqref{extra}. Assuming now that $u$ is a weak solution of \eqref{boundary-eqn-r} with $\lambda >0$. Let $M >0$ sufficiently large and will be determined. Let $u_M = u/M, \F_M = \F/M$ and $\lambda' = \lambda M$. We note that $u_M$ is a weak solution of
\begin{equation} \label{u-M-eqn}
\left\{
\begin{array}{cccl}
\text{div}[\A_{\lambda'} (x, u_M, \nabla u_M)] & = &\text{div}[\F_M], & \quad \text{in} \quad Q_{2R}^+, \\
u_M & = & 0, & \quad \text{on} \quad T_{2R} \times \Gamma_{2R}.
\end{array}
\right.
\end{equation}
Let us denote
\[
E_M = \Big\{Q_R^+: \M_{Q_{2R}^+}(|\nabla u_M|^2) > N \Big\}.
\]
and
\begin{equation} \label{K-zero}
K_0 = \left( \fint_{Q_{2R}^+ } |\nabla u|^2dxdt \right)^{1/2}.
\end{equation}
We claim that we can choose $M  = C K_0$ with some sufficiently large constant $C$ depending only on $\Lambda, M_0, M_1$, $M_2, p,  n$ such that
\begin{equation} \label{M-density}
\omega(E_M) \leq \epsilon \omega(Q_{R}(z)), \quad \forall \ z \in \overline{Q}_R^+.
\end{equation}
If this holds, we can apply \eqref{ineq-extra-cond-2} for $u_M$ which is a weak solution of \eqref{u-M-eqn} to obtain
\[
\norm{\nabla u_M}_{L^p(Q_{R}^+, \omega)} \leq C \left[ \norm{\F_M}_{L^p(Q_{2R}^+, \omega)}  +  \omega(Q_{R})^{1/p} \max\{(\lambda M)^{-1},1\}\right].
\]
Then, by multiplying this equality with $M$, we obtain
\[
\norm{\nabla u}_{L^p(Q_R^+, \omega)} \leq C \left[ \norm{\F}_{L^p(Q_{2R}^+, \omega)}  + \omega(Q_{R})^{1/p}\max\{\lambda^{-1}, K_0\}\right]
\]
which is \eqref{W-1p-flat} as desired. Therefore, the proof is therefore complete if we can prove \eqref{M-density}.  To this end, using the doubling property of $\omega \in A_q$, Lemma \ref{doubling}, we see that for every $z \in Q_1^+$,
\[
\begin{split}
\frac{\omega\big(E_M\big)}{\omega(Q_R(z))} & = \frac{\omega\big(E_M\big)}{\omega(Q_{2R})} \frac{\omega(Q_{2R})}{ \omega(Q_R(z))}  \leq   [\omega]_{A_q}  \frac{\omega\big(E_M\big)}{\omega(Q_{2R})} \left(\frac{|Q_{2R}|}{|Q_R(z)|}\right)^q = M_2(2)^{nq}   \frac{\omega\big(E_M\big)}{\omega(Q_{2R})}.
\end{split}
\]
Then, using Lemma \ref{doubling} again, we can find $\beta_0 = \beta_0 (M_2, n) >0$ such that
\[
\frac{\omega\big(E_2\big)}{\omega(Q_R(z))} \leq C(M_2, q, n) \left(\frac{|E|}{|Q_{2R}|} \right) ^{\beta_0}.
\]
On the other hand, by the weak type (1,1) estimate in Lemma \ref{Hardy-Max-p-p}, we see that
\[
\begin{split}
\big|E_M\big|  = \big|\big\{Q_R^+: \M_{Q_{2R}^+} (|\nabla u |^2) > N M^2 \big\}\big|  &  \leq \frac{C(n)}{N M^2} \int_{Q_{2R}^+}  |\nabla u|^2  dxdt \\
& =  \frac{C |Q_{2R}|}{M^2} \fint_{Q_{2R}^+} |\nabla u|^2  dxdt. 
\end{split}
\]
Hence, combining the last two estimates, we can find $C_*(M_2, q, n)  >1$ such that
\[
\begin{split}
\frac{\omega\big(E_M\big)}{\omega(Q_R(z))} & \leq C_*(M_2, q, n)  \left( \frac{1}{M^2}\fint_{Q_{2R}^+} |\nabla u|^2  dxdt   \right)^{\beta_0}, \quad \forall \ z \in Q_R^+.
\end{split}
\]
Then, by taking 
\begin{equation} \label{kappa-def}
M = \left[ \left( \fint_{Q_{2R}^+} |\nabla u|^2 dxdt \right)^{1/2} +1\right]\left( C_*/\e \right) ^{1/(2\beta_0)} = K_0 \left( C_*/\e \right) ^{1/(2\beta_0)} >1,
\end{equation} 
we obtain
\[
\omega(E_M) = \omega\big(\big\{Q_R^+ : \M_{Q_{2R}^+} (|\nabla u_M|^2) > N  \big\}\big) \leq \e \omega(Q_R(z)), \quad \forall \ z \in Q_R^+.
\]
This proves \eqref{M-density} and completes the proof.
\end{proof}
\section{Proof of global weighted $W^{1,p}$-regularity estimates} \label{main-theorm-proof}
\begin{proof}[\textbf{\textup{Proof of Theorem \ref{W-1p-reg}}}] The proof is standard once the two local interior and boundary regularity theorems are established, i.e. Theorem \ref{interior-thm} and Theorem \ref{W-1p-reg-b}. We mainly using partition of unity on $\Omega$ and then flatten the boundary $\partial \Omega$. Observe that the process of flattening the boundary $\partial \Omega$ will not significant change the $\textup{BMO}^\#$ of the coefficients as our domain is $C^1$. For details, one can find it in \cite{LTT}. By doing this, we then obtain 
our desired estimate
\[
\begin{split}
 \int_{\Omega \times(\bar{t}, T)}   |\nabla u|^p \omega(x,t) dx dt &  \leq C \left[ \int_{\Omega \times(0, T) } | \F|^p \omega(x, t) dx dt \right.\\
& \quad + \left. \omega(\Omega \times(0, T) ) \left\{ \fint_{\Omega \times(0, T)} |\nabla u|^2 dx dt\right)^{p/2}  +1 \right\}.
\end{split}
\]
For details, one can find it in \cite{LTT}. 
\end{proof}
\begin{remark}  \label{remark-2} Though, the initial data $u_0$ is defined in \eqref{main-eqn}, we do not require any condition on $u_0$ in this paper. Moreover, $\A(x, t, u, \nabla u)$ is not required to be continuous in $u$. Also, note that the estimates \eqref{CZ-est} and \eqref{W-1p-flat} are invariant under the dilation $u(x,t) \rightarrow u_{s}(x,t) : = u(sx, s^2t)/s$ and this is an intrinsic property in the Calder\'on-Zygmund theory.
\end{remark}
\noindent
\textbf{Acknowledgement.}  T. Phan's research is supported by the Simons Foundation, grant \#~354889.

\end{document}